\documentclass[12pt,a4paper]{article}

\pdfoutput=1
\usepackage{mystyle}
\usepackage{url}
\usepackage{cleveref}

\usepackage{amssymb}
\usepackage{amsthm}

\usepackage{bbm}

\def\R{\mathbb{R}}

\def\B{\mathcal{B}}

\def\conv{\mathrm{conv}}
\def\dim{\mathrm{dim}}

\def\ran{\mathrm{ran}}
\def\ker{\mathrm{ker}}
\def\lin{\mathrm{lin}}

\def\M{\mathcal{M}}

\def\one{\mathbbm{1}}

\urldef\urlpos\url{http://kyb.tuebingen.mpg.de/fileadmin/user_upload/files/publications/attachments/ilas_%5b0%5d.pdf}

\newtheorem{theorem}{Theorem}

\newtheorem{corollary}{Corollary}

\newtheorem{proposition}{Proposition}
\newtheorem{remark}{Remark}

\newcommand{\yo}[1]{{\color{cyan} \textbf{Yoh:}#1}}

\begin{document}

\title{On Representer Theorems and Convex Regularization}

\author{Claire Boyer$^1$, Antonin Chambolle$^2$, Yohann De Castro$^{3,4}$, Vincent Duval$^{4,5}$,  Fr\'ed\'eric de Gournay$^6$ and Pierre Weiss$^7$.\\
\footnotesize $^1$ LPSM, Sorbonne Universit\'e, ENS Paris
\footnotesize $^2$ CMAP, CNRS, Polytechnique, Palaiseau \\
\footnotesize $^3$ LMO, Universit\'e Paris-Sud, Orsay 
\footnotesize $^4$ MOKAPLAN, INRIA, Paris \\
\footnotesize $^5$ Ceremade, Universit\'e Paris Dauphine
\footnotesize $^6$ INSA, Universit\'e de Toulouse \\
\footnotesize $^7$ ITAV, CNRS, Universit\'e de Toulouse
}

\maketitle

\begin{abstract}
We establish a general principle which states that regularizing an inverse problem with a convex function yields solutions which are convex combinations of a small number of \textit{atoms}. These atoms are identified with the extreme points and elements of the extreme rays of the regularizer level sets. An extension to a broader class of quasi-convex regularizers is also discussed. As a side result, we cha\-rac\-terize the minimizers of the total gradient variation, which was still an unresolved problem.
\end{abstract}


\section{Introduction}

Let $\vecgal$ denote a {real} vector space.
Let $\Phi : \vecgal\to \RR^m$ be a linear mapping called \emph{sensing operator} and $u\in \vecgal$ denote a signal. 
The main results in this paper describe the structural properties of certain solutions of the following problem:
\begin{equation}\label{eq::mainproblem}
 \inf_{u \in \vecgal} f(\Phi u) + \reg(u), 
\end{equation}
where $R:\vecgal\to \R\cup\{+\infty\}$ is a convex function called \emph{regularizer} and $f$ is an arbitrary convex or non-convex function called \emph{data fitting term}. 
{
In many applications, one looks for ``sparse solutions'' that are linear sums of a few atoms. This article investigates the theoretical legitimacy of this usage. 
}

\paragraph{Representer theorems and Tikhonov regularization}

The name \emph{representer theorem} comes from the field of machine learning \cite{Scholkopf2002Learning}. 
To provide a first concrete example\footnote{Here, we follow the presentation of \cite{gupta2018continuous}.}, assume that $\Phi \in \RR^{m\times n}$ is a finite dimensional measurement operator and $L\in \R^{p\times n}$ is a linear transform. Solving an inverse problem using Tikhonov regularization amounts to finding the mini\-mi\-zers of
\begin{equation}\label{eq:tikhonov}
 \min_{u\in \RR^m} \frac{1}{2}\|\Phi u - y\|_2^2 + \frac{1}{2}\|Lu\|_2^2.
\end{equation}
Provided that $\ker\Phi\cap \ker L=\{0\}$, it is possible to show that, whatever the data $y$ is, solutions are always of the form
\begin{equation}\label{eq:firstrepresentertheorem}
 u^\star = \sum_{i=1}^m \alpha_i \psi_i + u_K,
\end{equation}
where $u_K\in \ker(L)$ and $\psi_i=(\Phi^T\Phi + L^TL)^{-1}(\phi_i)$, where $\phi_i^T\in \RR^n$ is the $i$-th row of~$\Phi$. This result characterizes structural properties of the minimi\-zers without actually needing to solve the problem. In addition, when $\vecgal$ is an infinite dimensional Hilbert space, Equation \eqref{eq:firstrepresentertheorem} sometimes allows to compute exact solutions, by simply solving a finite dimensional linear system. This is a critical observation that explains the practical success of kernel methods and radial basis functions \cite{Wendland2005Scattered}.

\paragraph{Representer theorems and convex regularization}

The Tikhonov regularization \eqref{eq:tikhonov} is a powerful tool when the number $m$ of observations is large and the operator $\Phi$ is not too ill-conditioned. However, recent results in the fields of compressed sensing \cite{Donoho2006Compressed}, matrix completion \cite{candes2009exact} or super-resolution \cite{Tang2013Compressed,candes2014towards} - to name a few - suggest that much better results may be obtained in general, by using convex regularizers, with level sets containing singularities. 
Popular examples of regularizers in the finite dimensional setting include the indicator of the nonnegative orthant \cite{donoho2005sparse}, the $\ell^1$-norm~\cite{Donoho2006Compressed} or its composition with a linear operator \cite{rudin1992nonlinear} and the nuclear norm \cite{candes2009exact}.
Those results were nicely unified in \cite{chandrasekaran2012convex} and one of the critical arguments behind all these techniques is a representer theorem of type \eqref{eq:firstrepresentertheorem}. 
In most situations however, this argument is only implicit. The main objective of this paper is to state a generalization of \eqref{eq:firstrepresentertheorem} to arbitrary convex functions $R$. 
It covers all the aforementioned problems, but also new ones for problems formulated over the space of measures.

To the best of our knowledge, the name ``{\it representer theorem}'' is new in the field of convex regularization and its first mention is due to Unser, Fageot and Ward in~\cite{unser2017splines}. Describing the solutions of~\eqref{eq::mainproblem} is however an old problem which has been studied since at least the 1940's in the case of Radon measure recovery.

\paragraph{Total variation regularization of Radon measures}
A typical example of inverse problem in the space of measures is
\begin{equation}\label{eq::beurling}
  \min_{\mu\in \Mm(\Omega)} |\mu|(\Omega)  \quad \mbox{s.t.}\quad  \Phi \mu=y 
\end{equation}
where $\Omega\subseteq \RR^N$, $\Mm(\Omega)$ denotes the space of Radon measures,  $|\mu|(\Omega)$ is the total variation of the measure $\mu$ (see Section~\ref{sec:app}) and $\Phi \mu$ is a vector of \textit{generalized moments}, \ie{} $\Phi \mu=\left(\int_\Omega \varphi_i(x)\d\mu(x)\right)_{1\leq i\leq m}$ where $\{\varphi_i\}_{1\leq i\leq m}$ is a family of continuous functions (which ``vanish at infinity'' if $\Omega$ is not compact). 

Problems of the form~\eqref{eq::beurling} have received considerable attention since the pioneering works of Beurling~\cite{Beurling1938} and Krein \cite{Krein1938}, sometimes under the name \textit{L-moment problem} (see the monograph~\cite{krein_markov_1977}). To the best of our knowledge, the first ``representer theorem'' for problems of the form~\eqref{eq::mainproblem} is given for~\eqref{eq::beurling} by \zuho{}~\cite{Zuhovickii1948} (see~\cite[Th. 3]{Zuhovickii1962} for an English version). It essentially states that
\begin{equation}\label{statement1}
  \mbox{\textit{There exists a solution to~\eqref{eq::beurling} of the form $\displaystyle\sum_{i=1}^{r}a_i \delta_{x_i}$, with $r\leq m$.}}
\end{equation}

{A more precise result was given by Fisher and Jerome in~\cite{fisher_spline_1975}. When considering the problem \eqref{eq::beurling}, and for a bounded domain $\Omega$, the result reads as follows:
\begin{align}\label{statement2}
\begin{split}
  &\mbox{\textit{The extreme points of the solution set to~\eqref{eq::beurling} are of the form}} \\
  &\qquad\qquad \qquad\qquad\qquad \displaystyle\sum_{i=1}^{r}a_i \delta_{x_i}, \mbox{  \textit{with}  } r\leq m.
\end{split}
\end{align}
Incidentally, the Fisher-Jerome theorem considers more general problems of the form:
\begin{equation}\label{eq::spline}
  \min_{u\in \vecgal} |Lu|(\Omega)  \quad \mbox{s.t.}\quad  Lu\in \Mm(\Omega) \qandq \Phi u =y,
\end{equation}
where $\vecgal\subseteq \mathcal{D}'(\Omega)$ is a suitably defined Banach space of distributions, $L:\mathcal{D}'(\Omega)\rightarrow \mathcal{D}'(\Omega)$ maps $\vecgal$ onto $\Mm(\Omega)$ and $\Phi:\vecgal\to \RR^m$ is a continuous linear operator. We refer to Section~\ref{sec:app} for precise assumptions. Let us mention that the initial results by Fisher-Jerome were extended to a significantly more general setting in \cite{unser2017splines}.}

It is important to note that the Fisher-Jerome theorem~\cite{fisher_spline_1975} provides a much finer description of the solution set than \zuho's result \cite{Zuhovickii1948}. Indeed, the well-known Krein-Milman theorem states that, if $\vecgal$ is endowed with the topology of a locally convex Hausdorff vector space and $\cvx\subset \vecgal$ is compact convex, then $\cvx$ is the closed convex hull of its extreme points, 
\begin{equation}
  \label{eq:krein}
  \cl \conv \left(\ext(\cvx)\right)=\cvx.
\end{equation}
In other words, the solutions described by the Fisher-Jerome theorem are sufficient to recover \emph{the whole set of solutions}.
Let us mention that the Krein-Milman theorem was extended by Klee~\cite{klee_extremal_1957} to unbounded sets: if $C$ is locally compact, closed, convex, and contains no line, then
\begin{equation}
  \label{eq:klee}
  \cl \conv \left(\ext(\cvx)\cup\rext(\cvx)\right)=\cvx,
\end{equation}
where $\rext(\cvx)$ denotes the union of the extreme rays of $\cvx$ (see Section~\ref{sec:notations} below).

\paragraph{``Representer theorems'' for convex sets}
As the Dirac masses are the extreme points of the total variation unit ball, each of the above-mentioned ``representer theorems'' for inverse problems actually reflect some phenomenon in the geometry of convex sets. 
In that regard, the celebrated Minkowski-Carath\'eodory theorem~\cite[Th.~III.2.3.4]{hiriart-urruty_convex_1993} is fundamental: any point of a compact convex set in an $m$-dimensional space is a convex combination of (at most) $m+1$ of its extreme points. 
In~\cite[Th.~(3)]{klee_theorem_1963}, Klee removed the boundedness assumption and obtained the following extension: any point of a closed convex set in an $m$-dimensional space is a convex combination of (at most) $m+1$ extreme points, or $m$ points, each an extreme point or a point in an extreme ray.

One purpose of the present paper is to point out the connection between the Fisher-Jerome theorem and a lesser known theorem by Dubins~\cite{dubins1962extreme} (see also~\cite[Exercise II.7.3.f]{bourbaki_espaces_2007}): 

\medskip

\begin{center}
\noindent\emph{The extreme points of the intersection of $\cvx$ with an affine space of codimension $m$ are convex combination of $($\!at most$)$\footnote{In the rest of the paper, we omit the mention ``at most'', with the convention that some points may be chosen identical.} $m+1$ extreme points of~$\cvx$}, 
\end{center}

\medskip

\noindent
provided $\cvx$ is linearly bounded and linearly closed (see Section~\ref{sec:notations}). That theorem was extended by Klee~\cite{klee_theorem_1963} to deal with the unbounded case. 

Although the connection with the Fisher-Jerome is striking, Dubins' theorem actually provides one extreme point too many. In the case of~\eqref{eq::beurling}, it would yield two Dirac masses for one linear measurement. We provide in this paper a refined analysis of the case of variational problems, which ensures at most $m$ extreme points.

\paragraph{Contributions}

The main results of this paper yield a description of some solutions to \eqref{eq::mainproblem} of the following form:
\begin{equation*}
 u^\star = \sum_{i=1}^r \alpha_i \psi_i + u_K,
\end{equation*}
where $r\leq m$, the atoms $\psi_i$ are identified with some extreme points (or points in extreme rays) of the regularizer level sets, and $u_K$ is an element of the so-called constancy space of $R$, \ie{} the set of directions along which $\reg$ is invariant. The results take the form \eqref{statement1}, when $f$ is an arbitrary function and the form \eqref{statement2} when it is convex. 
We provide tight bounds on the number of atoms $r$ that depend on the geometry of the level sets and on the link between the constancy space of $R$ and the measurement operator $\Phi$.

Our general theorems then allow us to revisit many results of the literature (linear programming, semi-definite programming, nonnegative constraints, nuclear norm, analysis priors), yielding simple and accurate descriptions of the minimizers. Our analysis also allows us to characterize the solutions of a resisting problem: we provide a representation theorem for the minimizers of the total gradient variation \cite{rudin1992nonlinear} as sums of indicators of simple sets. This provides a simple explanation to the staircaising effect when only a few measurements are used.

{Let us mention that, shortly after this work was posted on arXiv, similar results appeared, with somewhat different proofs, in a paper by Bredies and Carioni~\cite{bredies_sparsity_2018}.}

\section{Notation and Preliminaries}\label{sec:notations}
Throughout the paper, unless otherwise specified, $\vecgal$ denotes a finite or infinite dimensional real vector space and $\cvx\subseteq \vecgal$ is a convex set. Given two distinct points $x$ and $y$ in $\vecgal$, we let $\oi{x}{y}=\enscond{tx+(1-t)y}{0<t<1}$ and $[x,y]=\{tx+(1-t)y : $ $0\leq t\leq 1\}$ denote the open and closed segments joining $x$ to $y$. We recall the following definitions, and we refer to~\cite{dubins1962extreme,klee_extremal_1957} for more details. 



\paragraph{Lines, rays, and linearly closed sets} 
A \textit{line} is an affine subspace of $\vecgal$ with dimension $1$. An open half-line, \ie{} a set of the form $\rho=\{ p+tv : $ $t>0 \}$, where $p,v\in \vecgal$, $v\neq 0$, is called a \emph{ray} (through $p$). 
We say that the set $\cvx$ is \emph{linearly closed} (resp. linearly bounded) if the intersection of $\cvx$ and a line of $\vecgal$ is closed (resp. bounded) for the natural topology of the line. If $\vecgal$ is a topological vector space and $\cvx$ is closed for the corresponding topology, then $\cvx$ is linearly closed.

If $\cvx$ is linearly closed and contains some ray $\rho=p+\RR_+^*v$, it also contains the endpoint $p$ as well as the rays $q+\R_+v$ for all $q\in \cvx$. Therefore, if $\cvx$ contains a ray (resp. line), it recesses in the corresponding direction.

\paragraph{Recession cone and lineality space} The set of all $v\in\vecgal$ such that $\cvx+\RR_+^*v\subseteq \cvx$ is a convex cone called the \emph{recession cone of $\cvx$}, which we denote by $\rec{\cvx}$. If $\cvx$ is linearly closed then so is $\rec{\cvx}$, and  $\rec{\cvx}$ is the union of $0$ and all the vectors $v$ which direct the rays of $\cvx$. In particular, $\cvx$ contains a line if and only the vector space 
\begin{equation}
\label{eq::lineality space}
\lin(\cvx)\eqdef \rec{\cvx}\cap (-\rec{\cvx})
\end{equation}
is non trivial. The vector space $\lin(\cvx)$ 
 is called the \emph{lineality space} of $\cvx$. It corresponds to the largest vector space of invariant directions for $\cvx$. 
 
If $\vecgal$ is finite dimensional and $\cvx$ is closed, the recession cone coincides with the \emph{asymptotic cone}.

\paragraph{Extreme points, extremal rays, faces}
An \textit{extreme point} of $\cvx$ is a point $p\in \cvx$ such that $\cvx\setminus \{p\}$ is convex. 
An \textit{extremal ray} of $\cvx$ is a ray $\rho \in \cvx$ such that if $x,y\in \cvx$ and $\oi{x}{y}$ intersects $\rho$, then $\oi{x}{y}\subset \rho$. If $\cvx$ contains the endpoint~$p$ of $\rho$ (\eg{} if $\cvx$ is linearly closed), this is equivalent to $p$ being an extreme point of $\cvx$ and $\cvx\setminus \rho$ being convex.

  Following~\cite{dubins1962extreme,klee_theorem_1963}, if $p\in \cvx$, the smallest face of $\cvx$  which contains $p$ is the union of $\{p\}$ and all the open segments in $\cvx$ which have $p$ as an inner point. We denote it by $\face{p}{\cvx}$. The (co)dimension of $\face{p}{\cvx}$ is defined as the (co)dimension of its affine hull. The collection of all elementary faces, $\{\face{p}{\cvx}\}_{p\in\cvx}$, is a partition of $\cvx$.
Extreme points correspond to the zero-dimensional faces of $\cvx$, while extreme rays are (generally a strict subcollection of the) one-dimensional faces. 

\paragraph{Quotient by lines}
As noted above, if $\cvx$ is linearly closed, it contains a line if and only if the vector space $\lin(\cvx)$ defined in \eqref{eq::lineality space}
is nontrivial.  In that case, letting $W$ denote some linear supplement to $\lin(\cvx)$ 
, we may write 
\begin{equation}
\label{eq::quotientbyline}
\cvx= \qcvx+\lin(\cvx), \textrm{ with } \qcvx\eqdef \cvx\cap W
\end{equation}
and the corresponding decomposition is unique (\ie{} any element of $\cvx$ can be decomposed in a unique way as the sum of an element of $\qcvx$ and $\lin(\cvx)$). The convex set $\qcvx$ (isomorphic to the projection of $\cvx$ onto the quotient space $\vecgal/\lin(\cvx)$) is then linearly closed, and the decomposition of $\cvx$ in elementary faces is exactly given by the partition $\{\face{p}{\qcvx}+\lin(\cvx)\}_{p\in\qcvx}$, where $\face{p}{\qcvx}$ is the smallest face of $p$ in $\qcvx$. 

One may check that $\qcvx$ contains no line, as its recession cone $\rec{\qcvx}$, the projection of $\rec{\cvx}$ onto $W$ parallel to $\lin(\cvx)$, is a \emph{salient} convex cone.  

%



\section{Abstract representer theorems}
\label{sec:abstract}
\subsection{Main result}
Our main result describes the facial structure of the solution set to
  \begin{equation}\label{eq:thminreg}
    \min_{u\in\vecgal} \reg(u) \quad \mbox{s.t.}\quad \Phi u=y, \tag{$\Pp$}
\end{equation}
where $y\in \RR^m$, $\Phi:\vecgal \rightarrow \RR^m$ is a linear operator,  and $m\leq \dim \vecgal$, $m<+\infty$. 
In the following, let $\minval$ denote the \emph{optimal value} of~\eqref{eq:thminreg}, $\sol$ denote its \emph{solution set}, and~$\minset$ denote the \emph{corresponding level set} of $\reg$, 
\begin{align}
  \label{eq:defminset}
\minset\eqdef \enscond{u\in \vecgal}{\reg(u)\leq \minval}.
\end{align}

\begin{theorem}
\label{thm:first}
Let $\reg:\vecgal \to \RR\cup\{+\infty\}$ be a convex function. Assume that $\inf_\vecgal \reg< t^\star< +\infty$, that $\sol$ is nonempty and that the convex set $\minset$ is linearly closed and contains no line. 
Let $p\in\sol$ and let $j$ be the dimension of the face $\face{p}{\sol}$. Then $p$ belongs to a face of $\minset$ with dimension at most $m+j-1$. 

In particular, $p$ can be written as a convex combination of:
\begin{itemize}
  \item[$\circ$] $m+j$ extreme points of $\minset$,
  \item[$\circ$] or $m+j-1$ points of $\minset$, each an extreme point of $\minset$ or in an extreme ray of $\minset$.
\end{itemize}
%
Moreover, $\rec{\sol}=\rec{\minset}\cap \ker(\Phi)$ and therefore $\lin(\sol)=\lin(\minset)\cap \ker(\Phi)$.
\end{theorem}


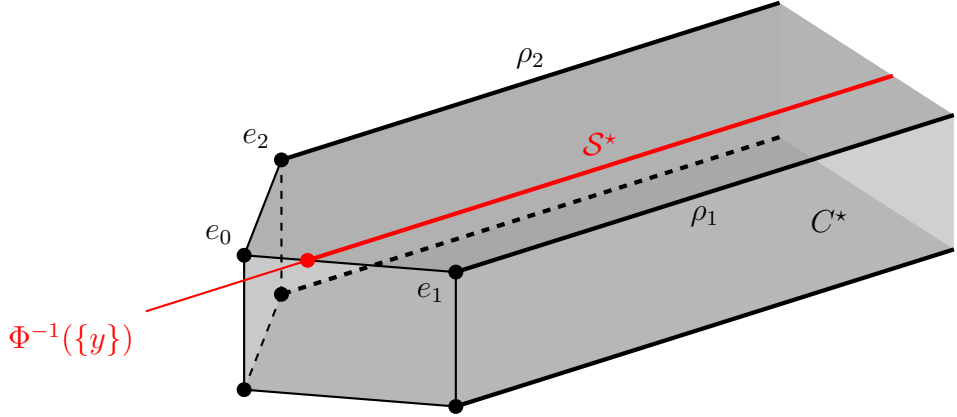
\begin{figure}
  \centering
  \tdplotsetmaincoords{63}{55}
    \begin{tikzpicture}[scale=2, tdplot_main_coords,axis/.style={->,dashed},thick]
\coordinate  (a0) at (0,0,1){};
\coordinate  (b0) at (0,0,0){};
\coordinate  (a1) at (1,1,1){};
\coordinate  (b1) at (1,1,0){};
\coordinate  (a2) at (1,5,1){};
\coordinate  (b2) at (1,5,0){};
\coordinate  (a3) at (-1,5,1){};
\coordinate  (b3) at (-1,5,0){};
\coordinate  (a4) at (-1,1,1){};
\coordinate  (b4) at (-1,1,0){};
\coordinate  (c0) at (0.3,-1,1){};
\coordinate  (c1) at (0.3,0.3,1){};
\coordinate  (c2) at (0.3,5,1){};

\coordinate  (d1) at (1,4,0.5){};

\fill[gray!80,opacity=0.7] (a0) -- (a1) -- (a2) -- (a3) -- (a4)  -- cycle; 
\fill[gray!80,opacity=0.5] (b0) -- (b1) -- (b2) -- (b3) -- (b4)  -- cycle; 
\fill[gray!80,opacity=0.4] (b0) -- (b1) -- (a1) -- (a0) -- cycle; 
\fill[gray!80,opacity=0.2] (b0) -- (b4) -- (a4) -- (a0) -- cycle; 
\fill[gray!80,opacity=0.4] (b1) -- (b2) -- (a2) -- (a1) -- cycle; 
\fill[gray!80,opacity=0.2] (b4) -- (b3) -- (a3) -- (a4) -- cycle; 
\fill[gray!80,opacity=0.1] (b2) -- (b3) -- (a3) -- (a2) -- cycle; 
\draw []       (a4)--(a0)--(a1);
\draw [ultra thick] (a1)--(a2) node[midway, below] {$\rho_1$}; 
\draw [ultra thick] (a3)--(a4) node[midway, above] {$\rho_2$};
\draw [dashed]      (b4)--(b0);
\draw [dashed,ultra thick]       (b3)--(b4);

\draw []       (b0)--(b1);
\draw [ultra thick] (b1)--(b2);

\draw []       (b1)--(a1);
\draw [dashed] (b4)--(a4);
\draw []       (b0)--(a0);

\draw [red,thick]       (c0)--(c1); %
\draw [red,ultra thick]       (c1)--(c2) node[midway, above] {$\sol$}; 
\draw[draw=red,fill=red] (c1) circle (0.1em);

\foreach \i in {0,1,4}
    {
      \draw[fill=black] (a\i) circle (0.1em);
      \draw[fill=black] (b\i) circle (0.1em);
    }
    \draw (a0) node[above left] {$e_0$};
      \draw (a1) node[below left] {$e_1$};
      \draw (a4) node[above left] {$e_2$};
      \draw[red] (c0) node[below left] {$\Phi^{-1}(\{y\})$};

      \draw[] (d1) node[] {$\minset$};

\end{tikzpicture}
    \caption{An illustration of \cref{thm:first} for $m=2$. The solution set $\sol=\minset\cap \Phi^{-1}(\{y\})$ is made of an extreme point and an extreme ray. The extreme point is a convex combination of~$\{e_0,e_1\}$. Depending on their position, the points in the ray are a convex combination of~$\{e_0,e_1,e_2\}$ or a pair of points, one in $\rho_1$ and the other in $\rho_2$.}
  \label{fig:thm} 
  \vspace{0.5cm}
\end{figure}

The proof of \cref{thm:first} is given in Section~\ref{sec:prooffirst}.
Before extending this theorem to a wider setting, let us formulate some remarks.

\begin{remark}[Extreme points and extreme rays of $\sol$]
  In particular ($j=0$), \emph{each extreme point} of $\sol$ is a convex combination of $m$ extreme points of~$\minset$, or a convex combination of $m-1$ points of~$\minset$, each an extreme point of~$\minset$ or in an extreme ray. Similarly $(j=1)$, \emph{each point on an extreme ray} of~$\sol$ is a convex combination of $m+1$ extreme points of~$\minset$, or a convex combination of $m$ points of~$\minset$, each an extreme point of $\minset$ or in an extreme ray. {Hence, provided the assumptions of Klee's theorem (see \eqref{eq:klee}) hold, Theorem \ref{thm:first} completely charaterizes the solution set.} An illustration is provided in \cref{fig:thm}.
\end{remark}

\begin{remark}[The hypothesis $\inf_\vecgal \reg< t^\star$\label{rem:infR}]
We have focused on the case $\minval>\inf_\vecgal \reg$ in the theorem since the case $\minval=\min \reg$ is easier. 
In that case, $M=\Phi^{-1}(\{y\})$ can be in arbitrary position (\ie{} not necessarily tangent) w.r.t.\ $\minset=\argmin \reg$, and one can only use the general Dubins-Klee theorem~\cite{dubins1962extreme,klee_theorem_1963} to describe their intersection. As a result the conclusions of \cref{thm:first} are slightly weakened, $p$ belongs to a face of $\minset$ with dimension $m+j$, and  one must add one more point in the convex combination (\eg{}, for $j=0$, each extreme point of $\sol$ is a convex combination of $m+1$ extreme points of $\minset$, or $m$ points\dots). 
\end{remark}

\begin{remark}[Gauge functions or semi-norms] A common practice in inverse problems is to consider positively homogeneous regularizers $\reg$, such as (semi)-norms or gauge functions of convex sets. In that case the extreme points of $\minset$ correspond, up to a rescaling, to the extreme points of $\{u\in\vecgal : \reg(u)\leq 1\}$. In several cases of interest, the extreme points of such convex sets are well understood, see Section~\ref{sec:app} for examples in Banach spaces or, for instance, the paper  \cite[Sec. 2.2]{chandrasekaran2012convex} for examples in finite dimensional spaces.
 \end{remark}
 
 \begin{remark}[Extension to semi-strictly quasi convex functions]
\Cref{thm:first} can be extended to the case where $R$ is a semi-strictly quasi-convex function. 
A function $R$ is said to be {\em semi-strictly quasi-convex} \cite{daniilidis2007some} if it is quasi-convex and if 
\[
  R(x)< R(y) \Longrightarrow R(\lambda x + (1-\lambda) y) < R(y) \quad \forall \lambda \in \oi{0}{1}.
\] 
In words, semi-strictly quasi-convex functions are functions that, when restricted to a line are successively decreasing, constant and increasing on their domain. In comparison, strictly quasi-convex functions are successively decreasing and increasing while quasi-convex functions successively non-increasing and non-decreasing. 

The set of semi-strictly quasi-convex functions is a subset of quasi-convex functions and it contains all convex and strictly quasi-convex functions. In the proof of \cref{thm:first}, only the semi-strictly quasi-convex property is required to ensure that \eqref{eq:for_cvx_like} holds.
 \end{remark}

 \begin{remark}[Topological properties\label{rem:lsc}] 
 The assumption that $\minset$ is linearly closed is fulfilled in most practical cases, since $\vecgal$ is usually endowed with the topology of a Banach (or locally convex) vector space and $\reg$ is assumed to be lower semi-continuous (so as to guarantee  the existence of a solution to~\eqref{eq:thminreg}). 
 {Note also that if $\reg$ is lower semi-continuous on any line (for the natural topology of the line), the set $\minset$ is linearly closed. 
 }
\end{remark}

%
%

\subsection{The case of level sets containing lines}
The reader might be intrigued by the assumption of \cref{thm:first} that $\minset$ contains no line, since in several applications the regularizer $\reg$ is invariant by the addition of, \eg, constant functions or low-degree polynomials (see Section~\ref{sec:app}). In that case, one is generally interested in the non-constant or non-polynomial part, and it is natural to consider a quotient problem for which the theorem applies. We describe below (see \cref{coro:lines}) how our result extends to the case where  $\minset$ contains lines. 

\begin{figure}[htbp]
  \centering
    \tdplotsetmaincoords{85}{12}
\begin{tikzpicture}[scale=1.4, tdplot_main_coords,axis/.style={->,dashed},thick]

\def\h{2.5} 
\def\xo{0} 
\def\yo{0} 
\def\zo{-1} 

\draw[thick,->,>=stealth] (\xo,\yo,\zo) -- (\xo+2.1*\h,\yo,\zo) node[anchor=north east]{};
\draw[thick,->,>=stealth] (\xo,\yo,\zo) -- (\xo,\yo+2.9*\h,\zo) node[anchor=north west]{};
\draw[thick,->,>=stealth] (\xo,\yo,\zo) -- (\xo,\yo,\zo+1+\h) node[anchor=south,left]{$\vecgal/\vecrec{}$};

\def\k{5} 
\coordinate  (a0) at (0,0,1){};
\coordinate  (b0) at (0,0,0){};
\coordinate  (a1) at (\k,0,1){};
\coordinate  (b1) at (\k,0,0){};
\coordinate  (a2) at (\k,\k,1){};
\coordinate  (b2) at (\k,\k,0){};
\coordinate  (a3) at (0,\k,1){};
\coordinate  (b3) at (0,\k,0){};

\fill[gray!80,opacity=0.7] (a0) -- (a1) -- (a2) -- (a3) --  cycle; 
\fill[gray!80,opacity=0.7] (b0) -- (b1) -- (b2) -- (b3) --  cycle; 
\fill[gray!80,opacity=0.5] (a0) -- (a1) -- (b1) -- (b0) --  cycle; 
\fill[gray!80,opacity=0.5] (b1) -- (b2) -- (a2) -- (a1) --  cycle; 
\fill[gray!80,opacity=0.5] (b2) -- (b3) -- (a3) -- (a2) --  cycle; 
\fill[gray!80,opacity=0.5] (b0) -- (b3) -- (a3) -- (a0) --  cycle; 

\draw [gray,dotted] (a0) -- (a1) -- (a2) -- (a3) --  cycle; 
\draw [gray,dotted] (b0) -- (b1) -- (b2) -- (b3) --  cycle; 
\draw [gray,dotted] (b1) -- (a1); 
\draw [gray,dotted] (b2) -- (a2); 
\draw [gray,dotted] (b3) -- (a3); 
\draw [gray,dotted] (b0) -- (a0); 
\draw[] (a1) node[below left] {$\minset$};

\def\a{0.5} 
\def\xshift{2.1} %

\coordinate  (c0) at (\xshift+\a*\h,0,\h){};
\coordinate  (c1) at (\xshift+\a*\h,\k,\h){};
\coordinate  (c2) at (\xshift-0.7*\a*\h,\k,-0.7*\h){};
\coordinate  (c3) at (\xshift-0.7*\a*\h,0,-0.7*\h){};

\coordinate  (c4) at (\xshift,0,0){};
\coordinate  (c5) at (\xshift,\k,0){};
\coordinate  (c6) at (\xshift+\a,0,1){};
\coordinate  (c7) at (\xshift+\a,\k,1){};

\draw[red,dotted] (c0) -- (c1) -- (c2) -- (c3) --  cycle; 
\fill[red!80,opacity=0.2] (c0) -- (c1) -- (c2) -- (c3) --  cycle; 

\draw [red,dotted, ultra thick]       (c5)--(c7); %
\draw [red,dotted, ultra thick]       (c4)--(c6); %
\fill [red!80,opacity=0.3]       (c4) -- (c6) -- (c7) --(c5) --cycle; %
\draw [red,ultra thick]       (c6)--(c7); %
\draw [red,ultra thick]       (c4)--(c5); %

\coordinate  (c8) at (\xshift+0.9*\a*\h,0.7*\k,0.7*\h){};
\draw[red] (c8) node [right] {$\Phi^{-1}(\{y\})$};
\draw[red] (c7) node[below right] {$\sol$};

\coordinate  (d0) at (\xo,\yo,0){};
\coordinate  (d1) at (\xo,\yo,1){};
\draw[teal,ultra thick] (d0) -- (d1) node[midway, left] {$\tilde\minset$}; 
\draw[draw=teal,fill=teal] (d0) circle (0.1em) ;
\draw[draw=teal,fill=teal] (d1)  circle (0.1em) ;

\draw[teal] (d0) node[left] {$q_1$};
\draw[teal] (d1) node[left] {$q_2$};

\coordinate  (k0) at (\xo,\yo,\zo){};
\coordinate  (k1) at (\xo+\k,\yo,\zo){};
\coordinate  (k2) at (\xo\k,\yo+\k,\zo){};
\coordinate  (k3) at (\xo,\yo+\k,\zo){};
\draw[blue,dotted] (k0) -- (k1) -- (k2) -- (k3) --  cycle; 
\fill[blue!80,opacity=0.2] (k0) -- (k1) -- (k2) -- (k3) --  cycle; 
\coordinate  (k4) at (\xo+0.7*\k,\yo+0.5*\k,\zo){};
\draw[blue] (k4) node[] {$\vecrec{}=\lin(\minset)$};

\draw[blue!50!red,dotted,thin] (\xshift+\a*\zo,\yo,\zo) -- (\xshift+\a*\zo,\yo+\k,\zo);
\end{tikzpicture}
    \caption{Taking the quotient by $\vecrec{}=\lin(\minset)$ yields a level set $\tilde\minset$ with no line. In this figure, to simplify the notation, we have omitted the isomorphism $\qiso$ (\ie{} in this figure $\minset$ shoud be replaced with $\qiso(\minset)$, and similarly for $\sol$ and $\Phi^{-1}(\{y\})$).}
  \label{fig:coro1b} 
\end{figure}


If $\minset$ is linearly closed and contains some line, it is translation-invariant in the corresponding direction. The collection of all such directions is the lineality space of $\minset$ (see \Cref{sec:notations}), we denote it by $\vecrec{}\eqdef \lin(\minset)$ (typically, if $\reg$ is the composition of a linear operator and a norm, $\vecrec{}$ is the \emph{kernel} of that linear operator). Let $\projv:\vecgal \rightarrow \vecgal/\vecrec{}$ be the canonical projection map. We recall that there exists a linear isomorphism $\qiso: \vecgal\rightarrow (\vecgal/\vecrec{})\times \vecrec{}$ such that the first component of $\qiso(p)$ is $\projv(p)$ for all $p\in \vecgal$.
We may now describe the equivalence classes (modulo $\vecrec{}$) of the solutions. 
\begin{corollary}
\label{coro:lines}
Let $\reg:\vecgal \to [-\infty,+\infty]$ be a convex function. Assume that $\inf_\vecgal \reg< t^\star< +\infty$, that $\sol$ is nonempty and that the convex set $\minset$ is linearly closed. Let $\vecrec{}\eqdef\lin(\minset)$ be the lineality space of $\minset$ and $d\eqdef \dim \Phi(\vecrec{})$. Let $p\in\sol$, let $\projv(p)$ denote its equivalence class, and let $j$ be the dimension of the face $\face{\projv(p)}{\projv(\sol)}$.

Then, $\projv(p)$ belongs to a face of $\projv(\minset)$ with dimension at most $m+j-d-1$.
In particular,
\begin{itemize}
  \item[$\circ$]  $\projv(p)$ is  a convex combination of $m+j-d$ extreme points of $\projv(\minset)$,
  \item[$\circ$] or $\projv(p)$ is a convex combination of $m+j-d-1$ points of $\projv(\minset)$, each an extreme point of $\projv(\minset)$ or in an extreme ray of $\projv(\minset)$.
\end{itemize}

As a result, letting $\qe_1,\ldots, \qe_r$ denote those extreme points (or points in extreme rays), 
\begin{equation}\label{eq:convcomb}
  p= \sum_{i=1}^r \theta_i \qiso^{-1}(\qe_i,0) + u_{\vecrec{}}, \qwhereq  \theta_i\geq 0,\  \sum_{i=1}^r \theta_i =1, \qandq u_{\vecrec{}}\in \vecrec{}.
\end{equation}
\end{corollary}
The proof of \cref{coro:lines} is given in \Cref{sec:proofconvex}.

One can have an explicit representation with elements of $E$ of a solution $p\in \sol $.
Indeed, let $W$ be some linear complement to $K=\lin(\minset)$. One may decompose $\minset = \tilde \minset + K$, where $\tilde{\minset}=\cvx\cap W$, and observe that  $\projv (\minset) $ and $\tilde \minset$ are isomorphic.
In this case, \cref{coro:lines} implies that
 $p$ can be written as the sum of one point in $\lin(\minset)$ and of a convex combination of:
\begin{itemize}
  \item[$\circ$] $m+j-d$ extreme points of $\tilde \minset$,
  \item[$\circ$] or $m+j-1-d$ points of $\tilde \minset$, each an extreme point of $\tilde \minset$ or in an extreme ray of $\tilde \minset$.
\end{itemize}

\subsection{Extensions to data fitting functions}
In this section, we discuss the extension of the above results to more general problems of the form
\begin{equation}
 \label{eq::data-fitting:convex}
 \inf_{u \in \vecgal} f(\Phi u)+\reg(u), \tag{$\Pp_f$}
\end{equation}
where $f:\RR^m\rightarrow \RR\cup\{+\infty\}$ is an arbitrary fidelity term.

\subsubsection{Convex data fitting term}
When $f$ is a \emph{convex} data fitting function $f$, we get the following result.

\begin{corollary} \label{cor:convex_fit}
  Assume that $f$ is convex and that the solution set $\sol_{f}$ of \eqref{eq::data-fitting:convex} is nonempty. Let $p\in \sol_{f}$ such that $\minset\eqdef \{u\in \vecgal, \reg(u)\leq \reg(p)\}$  is linearly closed, $\vecrec{}\eqdef\lin{\minset}$ and let $j$ be the dimension of the face $\face{\projv(p)}{\projv(\sol_{f})}$.

  If  $\inf_\vecgal \reg < \reg(p)$, then the conclusions of \cref{coro:lines} (or \cref{thm:first} if $\vecrec{}=\{0\}$) hold.

  If $\inf_\vecgal \reg = \reg(p)$, they hold with $1$ more dimension (see \cref{rem:infR}).
\end{corollary}
Let us recall that, in view of \cref{rem:lsc}, if $\vecgal$ is a topological vector space and $R$ is lower semi-continuous, $\minset$ is closed regardless of the choice $p$. 
\begin{proof}
Let $y=\Phi p$ and 
consider the following problem:
\begin{align}
\label{eq:singleton}
\tag{$\Pp_{\{y\}}$}
    \min_{u\in\vecgal } \reg(u) &\quad \mbox{s.t.}\quad 
    \Phi u = y.
    \end{align}
Let $\sol_{\{y\}}$ denote its solution set. It is a convex subset of $\sol_{f}$, with $p\in \sol_{\{y\}}$. Additionally, if $j$ is the dimension of  $\face{p}{\projv(\sol_f)}$ (resp. $\face{p}{\sol_f}$ if $\minset$ contains no line), then the face  $\face{p}{\projv(\sol_{\{y\}})}$  (resp. $\face{p}{\sol_{\{y\}}}$) has dimension at most $j$, since $\sol_{\{y\}}\subseteq \sol_{f}$.
It suffices to apply \cref{coro:lines} (resp. \cref{thm:first}) to obtain the result.
\end{proof}

\begin{remark}[The case of a strictly-convex function]
In the case when $f$ is strictly convex, it is known that $\Phi \sol_f$ is a singleton, which means that $\sol_{\{y\}}=\sol_f$.
\end{remark}

\begin{remark}[The case of quasi-convex functions]
The result actually holds whenever the solution set $\sol_f$ is convex. In particular, this property holds when $f$ is quasi-convex and $R$ is convex.
\end{remark}

\subsubsection{Non-convex function}
\label{subsec:nonconvex}
  In the general case, \ie{} when $f:\RR^m\rightarrow \RR\cup\{+\infty\}$ is an arbitrary function, it is difficult to describe the structure of the solution set. However, one may choose a solution $p_0$ (as before, provided it exists) and observe that it is also a solution to~\eqref{eq:thminreg} for $y\eqdef \Phi p_0$. Then, one may apply \cref{coro:lines}, but the difficult part is that the dimensions $j$ to consider are with respect to the solution set $\sol$ of the \emph{convex} problem~\eqref{eq:thminreg}. Nevertheless, if one is able to assert that the solution set $\sol$ has at least an extreme point $p$, then \cref{coro:lines} ensures that $p$ can be written in the form~\eqref{eq:convcomb}, where $r\leq m$ and the $\qe_i$'s are extreme points (or points in extreme rays) of $\minset$. Since $p$ must also be a solution to~\eqref{eq::data-fitting:convex}, one obtains that there exists a solution to~\eqref{eq::data-fitting:convex} of the form~\eqref{eq:convcomb}.

  \subsection{Ensuring the existence of extreme points}\label{sec:existextreme}
It is important to note that, in \cref{thm:first}, the existence a face of $\sol$ with dimension $j$ is not guaranteed (nor, for $j=0$, the existence of extreme points). The convex set $\sol$ might not even have any finite-dimensional face! 
For instance, let $\vecgal$ be the space of Lebesgue-integrable functions on $[0,1]$. If $\reg(u)=\int_0^1|u(x)|\d x$,  $\Phi:\vecgal\rightarrow \RR$ is defined by $u\mapsto \int_0^1 u(x)\d x$, and $y=1$, then 
  \begin{equation}
    \sol=\enscond{u\in\vecgal}{\int_0^1 |u(x)|\d x\leq 1 \qandq \int_0^1 u(x)\d x=1}.
  \end{equation}
It is possible to prove that such a set $\sol$ does not have any extreme point. As a consequence $\sol$ does not have any finite-dimensional face (otherwise an extreme point of the closure of a face would be an extreme point of $\sol$).

  However,  \cref{thm:first} (in fact the Dubins-Klee theorem~\cite{dubins1962extreme,klee_theorem_1963}) asserts that, \textit{if there is} a finite-dimensional face in $\sol$, then $\minset$ has indeed extreme points (and possibly extreme rays), and the convex combinations of such points generate the above-mentioned face.

  As a result, it is crucial to be able to assert a priori the existence of some finite-dimensional face for $\sol$, and this is where topological arguments come into play. If $\vecgal$ is endowed with the topology of a locally convex (Hausdorff)  vector space, the theorems~\cite[3.3 and 3.4]{klee_extremal_1957} which generalize the celebrated Krein-Milman theorem, state that 
  \emph{$\sol$ has an extreme point provided}
\begin{itemize}
    \item[$\circ$] $\sol$ is nonempty, convex,
    \item[$\circ$] $\sol$ contains no line,
    \item[$\circ$] and $\sol$ is closed, locally compact.
\end{itemize}
The last two conditions hold in particular if $\sol$ is compact. Moreover, as in \cref{coro:lines}, the second condition can be ensured by considering a suitable quotient map, provided it preserves the other topological properties (\eg{} if $\lin(\cvx)$ has a topological complement).

{ 
\begin{remark}
  Whereas local compactness is a very strong property for topological vector spaces (implying their finite-dimensionality, see~\cite[Th.~3,Ch.~1]{bourbaki_espaces_2007}), it is not so difficult to ensure the local compactness of $\sol$ in practice. Indeed, very often, even the existence of solutions is usually ensured using compactness arguments for a suitable weak or weak-* topology. The unbounded cases require more specific arguments, but let mention that there are examples of cones which are locally compact without being contained in any finite-dimensional vector space. In Section~\ref{sec:momprob} below, we discuss the example of the cone $\M^+(\Omega)$ of non-negative measures over a compact set for the weak-* topology. Another example of locally compact convex cone is 
\[
  \mathcal C=\Big\{ x\in\RR^\NN  \text{ such that } x_n \ge 0 \text{ and } \sum_{n\in\mathbb{N} } x_n \omega_n \le  \sum_{n\in\mathbb{N} } x_n <+\infty \Big\}\subseteq \ell^1(\NN)\,,
\]
for some non-decreasing positive sequence $(\omega_n)_n$ converging to $+\infty$ (note that $\omega_0<{1}$ for the cone to be non-empty). The cone $\mathcal{C}$ is locally compact for the strong topology. 
Indeed, consider the intersection $K$ of the cone $\mathcal C$ and the strong unit ball, namely $K = \{x=(x_n)_n\ : \ x_n\geq0,\ \sum_n \omega_n x_n \le  \sum_{n } x_n\leq 1\}$ and consider a sequence of elements of $K$ denoted~$(x^k)_k \subset K$. Using a diagonal argument, each $(x_n^k)_k$ converges to some $\bar{x}_n \geq0$. Furthermore, using that $\{n\,:\, w_n<1\}$ is finite and Fatou's lemma, it holds that
\begin{align*}
\sum_{n } \bar x_n&\leq \liminf_k \sum_{n } x^k_n\leq1\quad \mathrm{and}\\
\sum_{n\,:\, w_n\geq1 } (w_n-1)\bar x_n&\leq \liminf_k \sum_{n\,:\, w_n\geq1 } (w_n-1) x^k_n\\
&\leq \liminf_k \sum_{n\,:\, w_n<1 } (1-w_n) x^k_n=\sum_{n\,:\, w_n<1 } (1-w_n)\bar x_n
\end{align*}
and we deduce that $\bar x= (\bar{x}_n)_n\in K$. Furthermore, one has
\[
||x^k-\bar{x}||_1 = \sum_{n=0}^M |x_n^k-\bar{x}_n| + \sum_{n>M} |x_n^k-\bar{x}_n|
\]
and $M$ will be chosen later. Finally,
\begin{align*}
  \sum_{n>M} |x_n^k-\bar{x}_n| &\leq \sum_{n>M} (x_n^k+ \bar{x}_n) \leq  (1/\omega_M) \sum_{n>M} \omega_n (x_n^k+\bar{x}_n)
\leq  2/\omega_M
\end{align*}
since $\omega_n/\omega_M \geq 1$ for $n>M$. Therefore choosing $M$ large enough ensures that the second term $\sum_{n>M} |x_n^k-\bar{x}_n|$ is less than some $\varepsilon>0$. Choosing $k$ large enough leads to $||x^k-\bar{x}||_1\leq 2\varepsilon$.
\end{remark}
}


\section{Application to some popular regularizers}
\label{sec:app}

We now show that the extreme points and extreme rays of numerous convex regularizers can be described analytically, allowing to describe important analytical properties of the solutions of some popular problems. 
The list given below is far from being exhaustive, but it gives a taste of the diversity of applications targeted by our main results.


\subsection{Finite-dimensional examples}
We first consider examples where one has $\dim \vecgal<+\infty$. In that case, $\Phi$ is continuous, and since the considered regularizations~$\reg$ are lower semi-continuous, we deduce that
\begin{itemize}
  \item[$\circ$] the level set $\minset$ is closed,
  \item[$\circ$] the solution set $\sol=\minset\cap \Phi^{-1}(\{y\})$ is closed, and locally compact (even compact in most cases), hence it admits extreme points provided it contains no line.
\end{itemize}

\subsubsection{Nonnegativity constraints}
\label{sec:nonneg}

In a large number of applications, the signals to recover are known to be nonnegative. In that case, one may be interested in solving nonnegatively constrained problems of the form:
\begin{equation}\label{eq:nonnegative}
 \inf_{u\in \RR^n_+} f(\Phi u - y).
\end{equation}
An important instance of this class of problems is the nonnegative least squares \cite{lawson1995solving}, which finds it motivation in a large number of applications. Applying the results of \Cref{sec:abstract} to Problem \eqref{eq:nonnegative} yields the following result.
\begin{proposition}
 If the solution set of \eqref{eq:nonnegative} is nonempty, then it contains a solution which is $m$-sparse.
 In addition if $f$ is convex and the solution set is compact, then its extreme points are $m$-sparse.
\end{proposition}
Choosing $\reg$ as the characteristic function of ${\RR^n_+}$, the result simply stems from the fact that the extreme rays of the positive orthant are the half lines $\{\alpha e_i, \alpha \geq 0\}$, where $(e_i)_{1\leq i \leq n}$ denote the elements of the canonical basis. We have to consider $m$ atoms and not $m-1$ since $t^\star=\inf \reg=0$, see \cref{rem:infR}.

It may come as a surprise to some readers, since the usual way to promote sparsity consists in using $\ell^1$-norms. This type of result is one of the main ingredients of \cite{donoho2005sparse} which shows that the $\ell^1$-norm can sometimes be replaced by the indicator of the positive orthant when sparse \emph{positive} signals are looked after.

\subsubsection{Linear programming}

Let $\psi \in \R^n$ be a vector and $\Phi \in \RR^{m\times n}$ be a matrix and consider the following linear program in standard (or equational) form:
\begin{equation}\label{eq:linearprogramming}
 \inf_{\substack{u \in \RR^n_+ \\ \Phi u = y}} \langle \psi,u\rangle
\end{equation}

Applying Theorem \eqref{thm:first} to the problem \eqref{eq:linearprogramming}, we get the following well-known result (see e.g. \cite[Thm. 4.2.3]{matousek2007understanding}):
\begin{proposition}\label{prop:linprog}
Assume that the solution set of \eqref{eq:linearprogramming} is nonempty and compact.
Then, its extreme points are $m$-sparse, i.e. of the form:
\begin{equation}
 u = \sum_{i=1}^m \alpha_i e_i, \alpha_i\geq0,
\end{equation}
where $e_i$ denotes the $i$-th element of the canonical basis.
\end{proposition}
In the linear programming literature, solutions of this kind are called \emph{basic solutions}.
To prove the result, we can reformulate \eqref{eq:linearprogramming} as follows:
\begin{equation}\label{eq:linearprogramming2} 
 \inf_{\substack{(u,t) \in \RR^n_+\times \RR \\ \Phi u = y \\ \langle \psi,u\rangle=t}} t
\end{equation}
\begin{sloppypar}
Letting $R(u,t) = t + \iota_{\RR_+^n}(u)$, we get $\inf R = -\infty$. Hence, if a solution exists, we only need to analyze the extreme points and extreme rays of ${\minset = \{(x,t)\in \RR^n\times \RR, R(x,t) \le t^\star\}=\RR_+^n\times ]-\infty,t^\star]}$, where $t^\star$ denotes the optimal value. The extreme rays of this set (a shifted nonnegative orthant) are of the form $\{\alpha e_i, \alpha > 0 \}\times \{t^\star\}$ or $\{0\}\times]-\infty,t^\star[$. In addition $\minset$ possesses only one extreme point $(0,t^\star)$. Applying \cref{thm:first}, we get the desired result.
\end{sloppypar}
\subsubsection{$\ell^1$ analysis priors}

An important class of regularizers in the finite dimensional setting $\vecgal=\RR^n$ contains the functions of the form $R(u) = \| Lu \|_1$, where {$L$ is a linear operator from $\RR^{n}$ to $\RR^{p}$}. They are sometimes called \emph{analysis priors}, since the signal $u$ is ``analyzed'' through the operator $L$. Remarkable practical and theoretical results have been obtained using this prior in the fields of inverse problems and compressed sensing, even though many of its properties are -to the belief of the authors- still quite obscure. 

Since $R$ is one-homogeneous, it suffices to describe the extremality properties of the unit ball $C=\{u\in \RR^n, \| Lu \|_1\leq 1\}$ to use our theorems. The lineality space is simply equal to $\lin(C) = \ker(L)$. Let $K=\ker(L)$, $K^\perp$ denote the orthogonal complement of $K$ in $\RR^n$ and $L^+:\RR^n\to K^\perp$ denote the pseudo-inverse of $L$.  
We can decompose $C$ as $C= K + C_{K^\perp}$ with $C_{K^\perp}=C\cap K^\perp$.
Our ability to characterize the extreme points of $C_{K^\perp}$ depends on whether $L$ is surjective or not.
Indeed, we have 
\begin{equation}\label{eq:extL1}
\ext(C_{K^\perp})=L^+\left( \ext\left( \ran(L)\cap B_1^p\right) \right),
\end{equation}
where $B_1^p$ is the unit $\ell^1$-ball defined as
\begin{equation*}
B_1^p=\{z\in \RR^p, \|z\|_1\leq 1\}. 
\end{equation*}
Property \eqref{eq:extL1} simply stems from the fact that $C_{K^\perp}$ and $D=\ran(L)\cap B_1^p$ are in bijection through the operators $L$ and $L^+$.

\paragraph{The case of a surjective operator $L$}

When $L$ is surjective {(hence $p\leq n$)}, the problem becomes quite elementary.
\begin{proposition}\label{prop:extanalysis}
If $L$ is surjective, the extreme points $u$ of $C_{K^\perp}$ are $\ext(C_{K^\perp}) = (\pm L^+ e_i)_{1\leq i \leq p}$, where $e_i$ denotes the $i$-th element of the canonical basis. 
{Consider Problem \eqref{eq::data-fitting:convex} and assume that at least one solution exists. } 
Then Problem \eqref{eq::data-fitting:convex} has solutions of the form
 \begin{equation}\label{eq:extremeanalysisfinite}
  u^\star = \sum_{i\in I} \alpha_i L^+e_i + u_K,
 \end{equation}
 where $u_K\in \ker(L)$ and $I\subset \{1,\hdots, p\}$ is a set of cardinality $|I|\leq m - \dim(\Phi\ker(L))$.
\end{proposition}
{The proof of Proposition \ref{prop:extanalysis} follows from Corollary \ref{coro:lines} and Section \ref{subsec:nonconvex}, with $j=0$ and observing that $\projv$ is the orthogonal projection on $K^\perp$.}

\paragraph{The case of an arbitrary operator $L$}

When $L$ is not surjective the des\-crip\-tion of the extreme points $\ext\left(D\right)$ becomes untractable in general. A rough upper-bound on the number of extreme points can be obtained as follows. 
We assume that $L$ has full rank {$n$} and that $\ran(L)$ is in general position. The extreme points of $\ran(L)\cap B_1^p$ correspond to the intersections of some faces of the $\ell^1$-ball with a subspace of dimension {$n$}. In order for some $k$-face to intersect the subspace $\ran(L)$ on a singleton, $k$ should satisfy $n+k-p=0$, i.e. $k=p-n$. The $k$-faces of the $\ell^1$-ball contain $(k+1)$-sparse elements. The number of $(k+1)$-sparse supports in dimension $p$ is $\binom{p}{k+1}$. For a fixed support, the number of sign patterns is upper-bounded by $2^{k+1}$. Hence, the maximal number of extreme points satisfies $|\ext\left( \ran(L)\cap B_1^p\right)|\leq 2^{k+1}\binom{p}{k+1}$. This upper-bound is pessimistic since the subspace may not cross all extreme points, but it provides an idea of the combinatorial explosions that may happen in general. 

Notice that enumerating the number of faces of a polytopes is usually a hard problem. For instance, the Motzkin conjecture \cite{motzkin1957comonotone} which upper bounds the number of $k$ faces of a $d$ polytope with $z$ vertices was formulated in 1957 and solved by Mc Mullen \cite{mcmullen1970maximum} in 1970 only.


\subsubsection{Matrix examples}
\label{sec:matrix:examples}
In several applications, one deals with optimization problems in matrix spaces. The following regularizations/convex sets are commonly used. 

\paragraph{Semi-definite matrix constraint}
Similarly to \Cref{sec:nonneg}, one may consider  in $\RR^{n\times n}$ the following constrained problem
\begin{equation}\label{eq:sdpcone}
 \inf_{M\succeq 0} f(\Phi M - y),
\end{equation}
where $M\succeq 0$ means that $M$ must be symmetric positive semi-definite (p.s.d.). The extreme rays of the positive semi-definite cone $\minset$ are the p.s.d. matrices of rank $1$ (see for instance~\cite[Sec. 2.9.2.7]{dattorro_convex_2005}).
Hence, arguing as in \Cref{sec:nonneg}, we may deduce that if there exists a solution to~\eqref{eq:sdpcone}, there is also a solution which has rank (at most) $m$.

However, that conclusion is not optimal, as in that case a theorem by Barvinok~\cite[Th. 2.2]{barvinok_problems_1995} ensures that there exists a solution $M$ with
\begin{equation}\label{eq:ranksdp}
    \rank(M) \leq \frac{1}{2}\left(\sqrt{8m+1}-1\right).
  \end{equation}
  To understand the gap with Barvinok's result, let us note that the p.s.d.\ cone has a very special structure which makes the Minkowski-Carath\'eodory theorem (or its extension by Klee) too pessimistic. By~\cite[Sec. 2.9.2.3]{dattorro_convex_2005}, given $M\succeq 0$, the smallest face of the p.s.d. cone which contains $M$ (\ie{} the set of p.s.d matrices which have the same kernel) has dimension 
  \begin{equation}\label{eq:ranksdp2}
    d= \frac{1}{2}\rank(M)(\rank(M)+1).
  \end{equation}
  \begin{sloppypar}
  Equivalently, if the smallest face which contains $M$ has dimension $d$, then $\rank(M)=\frac{1}{2}\left(\sqrt{8d+1}-1\right)$, hence $M$ is a convex combination of ${\frac{1}{2}\left(\sqrt{8d+1}-1\right)}$ points in extreme rays, a value which is less than the value $d$ predicted by Klee's extension of Carath\'eodory's theorem.

  As a result, we recover Barvinok's result by noting that, as ensured by the first claim\footnote{or, more precisely, by its variant when $t^\star=\inf \reg$, see \cref{rem:infR}.} of \cref{thm:first}, any extreme point $M$ of the solution set belongs to a face of dimension $m$. Then, taking into account~\eqref{eq:ranksdp2} improves upon the second claim of \cref{thm:first}, and we immediately obtain~\eqref{eq:ranksdp}. 
  \end{sloppypar}


\paragraph{Semi-definite programming}

Semi-definite programs are problems of the form:
\begin{equation}\label{eq:semidefiniteprogramming}
 \inf_{\substack{M\succeq 0\\ \Phi(M) = y}} \langle A,M\rangle,
\end{equation}
where $A\in \RR^{n\times n}$ is a matrix and $\langle A,M\rangle\eqdef \mathrm{Tr}(A M)$. 
Arguing as in \cref{prop:linprog}, if the solution set of \eqref{eq:semidefiniteprogramming} is nonempty, our main result allows to state that its extreme points are matrices of rank $m$. In view of the above discussion, it is possible to refine this statement and show that~\eqref{eq:ranksdp} holds.

\paragraph{The nuclear norm}
The nuclear norm of a matrix $M\in \RR^{p\times n}$ is often denoted $\|M\|_*$ and defined as the sum of the singular values of $M$. It gained a considerable attention lately as a regularizer thanks to its applications in matrix completion \cite{candes2009exact} or blind inverse problems \cite{ahmed2014blind}. 
The geometry of the unit ball $\{M\in  \RR^{p\times n}, \|M\|_*\leq 1\}$ is well studied due to its central role in the field of semi-definite programming \cite{pataki2000geometry}. Its extreme points are the rank one matrices $M= uv^T$, with $\|u\|_2=\|v\|_2=1$.

Combining \cref{thm:first} with this result explains why regularizing pro\-blems over the space of matrices with the nuclear norm allows recovering \emph{rank-m solutions}.

\paragraph{The rank-sparsity ball}
The \emph{rank-sparsity} ball is the set $\{M\in  \RR^{m\times n}, \|M\|_* + \|M\|_1 \leq 1\}$, where $\|M\|_1$ is the $\ell^1$-norm of the entries of $M$. The corresponding regularization is sometimes used in order to favor sparse and low-rank matrices. 
The authors of \cite{drusvyatskiy2015extreme} have described the extreme points of this unit ball.
They have proved that \emph{the extreme points $M$ of the rank sparsity ball satisfy $\frac{r(r+1)}{2}-|I|\leq 1$}, where $|I|$ denotes the number of non-zero entries in $M$ and $r$ denotes its rank.
This result partly explains why using the rank-sparsity gauge promotes sparse and low rank solution. Let us outline that this effect might be better obtained using different strategies \cite{Richard2013Intersecting,richard2014tight}.


\paragraph{Bi-stochastic matrices}
A doubly stochastic matrix is a matrix with nonnegative rows and columns summing to one. The set of such matrices is called the Birkhoff polytope. The Birkhoff-von Neumann theorem states that its extreme points are the permutation matrices. 
We refer the interested reader to \cite{Fogel2015Convex} for an use of such matrices in DNA sequencing.

\subsection{Examples in infinite dimension}

In this section, we provide results in infinite dimensional spaces, which echoe the ones described in finite dimension.

\subsubsection{Problems formulated in Hilbert or Banach sequence spaces}

The case of Hilbert spaces (or countable sequences) can be treated within our formalism and all the examples presented previously have their na\-tu\-ral counterpart in this setting. {In the same vein, one can also treat Banach sequence spaces $\ell^p$ for $1\leq p\leq\infty$}. We do not reproduce the results here for space limitations. Let us however mention that two works treat this specific case with $\ell^1$ regularizers \cite{unser2016representer,adcock2016generalized}.

\subsubsection{Linear programming and the moment problem}\label{sec:momprob}
~

Let $\Omega$ be a {compact} metric space, $\mathcal{M}(\Omega)$ be the set of Radon measures on $\Omega$ and let $\mathcal{M}_+(\Omega)\subseteq \mathcal{M}(\Omega)$ be the cone of nonnegative measures on $\Omega$.
Let $\psi$ and $(\phi_i)_{1\leq i \leq m}$ denote a collection of continuous functions on $\Omega$. 
Now, let $\Phi:\mathcal{M}(\Omega)\to \RR^m$ be defined by $(\Phi \mu)_i = \langle \phi_i, \mu\rangle$, where $ \langle \phi_i, \mu\rangle\eqdef\int_\Omega \phi_i \d \mu$,  and consider the following linear program in standard form:
\begin{equation}\label{eq:linearprogramminginf}
  \inf_{\substack{\mu \in \mathcal{M}_+(\Omega) \\ \Phi \mu = y}} \langle \psi,u\rangle.
\end{equation}

Applying Theorem \eqref{thm:first} to the problem \eqref{eq:linearprogramminginf}, we get \cref{prop:linearprogramminginf} below. We do not provide a proof here, since it mimics very closely the one given for linear programming in finite dimension. The extreme rays of $\mathcal{M}_+(\Omega)$ can be described, arguing as in~\cite[Th.~15.9]{aliprantis_infinite_2006}, as the rays directed by the Dirac masses.
\begin{proposition}\label{prop:linearprogramminginf}
  Assume that the solution set \eqref{eq:linearprogramminginf} is {nonempty}.
Then, its extreme points are $m$-sparse, i.e. of the form:
\begin{equation}
 u = \sum_{i=1}^m \alpha_i \delta_{x_i}, x_i\in \Omega, \alpha_i\geq0.
\end{equation}
\end{proposition}
{
  To make sure that the above proposition is non-trivial, one may wish to ensure that the solution set $\sol$ has indeed extreme points, using arguments from Section~\ref{sec:existextreme}. It is straightforward that $\sol$ is convex and does not contain any line. Now, let us endow $\mathcal{M}(\Omega)$ with the weak-* topology (\ie{} the coarsest topology for which $\mu\mapsto \int_{\Omega}\eta\d \mu$ is continuous for every $\eta\in \Cder{}(\Omega)$). By lower semi-continuity, $\sol$ is closed. Moreover, $\sol$ is locally compact since the closed convex cone $\mathcal{M}_+$ is itself locally compact (take any $\mu\in \mathcal{M}_+(\Omega)$, its neighborhood $\enscond{\nu\in \mathcal{M}_+(\Omega)}{\nu(\Omega)\leq \mu(\Omega)+1}$ is compact in the weak-* topology).
}

{Proposition~\ref{prop:linearprogramminginf}} is well known, see e.g. \cite{shapiro2001duality}. Note that if we optimize the linear form $\langle \psi,u \rangle$ over the set of probability measures instead of the set of nonnegative measures, we get the so-called moment problem \cite{shohat1943problem} for which we can obtain a similar result.

\subsubsection{The total variation ball}

Let $\Omega$ denote an open subset of $\RR^d$ and $\M(\Omega)$ denote the set of Radon measures on $\Omega$.
The total variation ball $B_\M = \{u \in \M(\Omega), \|u\|_{\M(\Omega)}\leq 1\}$ plays a critical role for problems such as super-resolution \cite{candes2014towards,Tang2013Compressed,Duval2015Exact}. {It is compact for the weak-* topology and }its extreme points are the Dirac masses: $\ext(B_\M)=\{\pm \delta_x, x\in \Omega\}$. Hence total variation regularized problems of the form:
\begin{equation*}
 \inf_{u \in \M} f(\Phi u) + \|u\|_\M,
\end{equation*}
yield $m$-sparse solutions (under an existence assumption). A few variations around this central result were provided in \cite{fernandez2016super}.

\paragraph{Demixing of Sines and Spikes}
In \cite[Page 262]{fernandez2016super}, the author presents a regularization of the type
\[
\|\mu\|_{\mathcal{M}}+\eta\|v\|_1
\]
where $\eta>0$ is a tuning parameter, $\mu$ a complex measure and $v\in\mathbb C^n$ a sparse vector. Define $E$ as the set of $(\mu,v)$ where $\mu$ is a complex Radon measure on a domain $\Omega$ and $v\in\mathbb C^n$. Consider the unit ball
\[
B\eqdef\{(\mu,v)\in E\ :\ \|\mu\|_{\mathcal{M}}+\eta \|v\|_1\leq 1\}.
\]
Its extreme points are the points 
\begin{itemize}
\item 
$(a\delta_t,0)$ for all $t\in\Omega$ (and $\delta_t$ denotes the Dirac mass at point $t$) and all $a\in\mathbb C$ such that $|a|=1$, 
\item
$(0,ae_k)$ for all $k=1,\ldots,n$ and all $a\in\mathbb C$ such that $|a|=1/\eta$ and $e_k$ denotes the vector with $1$ at entry $k$ and $0$ otherwise. 
\end{itemize}

\paragraph{Group Total Variation: Point sources with a common support}
In \cite[Page 266]{fernandez2016super}, the author presents a regularization of the type
\[
\|\mu\|_{\mathcal{M}^n}:=\sup_{F:\Omega\to\mathbb C^n,\ \| F(t) \|_2\leq 1,\ t\in\Omega}\int_\Omega \langle F(t), \nu(t)\rangle \mathrm d|\mu|(t)
\]
were $F$ is continuous and vanishing at infinity, and $\mu$ is a vectorial Radon measure on~$\Omega$ such that $|\mu|$-a.e. $\mu=\nu \cdot |\mu|$ with~$\nu$ a measurable function from~$\Omega$ onto $\mathbb S^{n-1}$ the $n$-sphere and $|\mu|$ a positive finite measure on $\Omega$. Consider the unit ball
\[
B\eqdef\{\mu, \|\mu\|_{\mathcal{M}^n}\leq 1\}.
\]
Its extreme points are $a\delta_t$ for all $t\in\Omega$ (and $\delta_t$ denotes the Dirac mass at point $t$) and all $a\in\mathbb C^n$ such that $\|a\|_2=1$.

\subsubsection{Analysis priors in Banach spaces}
\label{sec:prior}

The analysis of extreme points of analysis priors in an infinite dimensional setting is more technical. 
Fisher and Jerome~\cite{fisher_spline_1975} proposed an inte\-res\-ting result, which can be seen as an extension of \eqref{eq:extremeanalysisfinite}. This result was recently revisited in \cite{unser2017splines} and \cite{flinth2017exact}. Below, we follow the presentation in \cite{flinth2017exact}.

{
Let $\Omega$ denote an open set in $\RR^d$. 
Let $\mathcal{D}'(\Omega)$ denote the set of distributions on~$\Omega$ and let $L:\mathcal{D}'(\Omega) \to \mathcal{D}'(\Omega)$ denote a linear operator with kernel $K=\ker(L)$.
We let $E=\{u\in \mathcal{D}'(\Omega), Lu\in \M(\Omega)\}$ and let $\|\cdot\|_K$ denote a semi-norm on $E$, which restricted to $K$ is a norm. We define a function space $\mathcal{B}(\Omega)$ as follows: 
$$\mathcal{B}(\Omega) = \{u\in E, \|Lu\|_{\M(\Omega)}+\|u\|_K<+\infty\}$$
 and equip it with the norm $\|u\|_{\B(\Omega)} = \|Lu\|_{\M(\Omega)} + \|u\|_K$.
 We assume that $L$ is surjective, i.e. $\M(\Omega) = L(\mathcal{B}(\Omega))$, and that $K$ has a topological complement (with respect to $\mathcal{B}(\Omega)$), which we denote by $K^\perp$. This setting encompasses all surjective Fredholm operators for instance.
Under the stated assumptions, we can define a pseudo-inverse $L^+$ of $L$ relative to $K^\perp$ \cite{beutler}.

The representer theorems in \cite{fisher_spline_1975,unser2017splines,flinth2017exact} can be obtained using \cref{thm:first} as exemplified below.
 
\begin{proposition}
Let $B=\{u \in \B(\Omega), \|Lu\|_{\M(\Omega)}\leq 1\}$. 
Then the extreme points of the set $C_{K^\perp}=B\cap K^\perp$ are of the form $\pm L^+\delta_{x}$, for $x\in \Omega$.

Let $f:\R^m\to \R\cup \{+\infty\}$ denote a convex function and define
\begin{equation*}
\sol = \argmin_{u \in \B(\Omega)} f(\Phi u) + \|Lu\|_{\M(\Omega)}.
\end{equation*}
Assume that $\sol$ is nonempty and does not contain $0$. Then the extreme points (if they exist) of $\pi_K(\sol)$ are of the form $u = \sum_{i=1}^m \alpha_i L^+ \delta_{x_i}$.
\end{proposition}
\begin{proof}
The proof mimics the finite dimensional case \eqref{eq:extremeanalysisfinite}.
First notice that $B=L^{-1}(B_\M)$, where $L^{-1}(\{\mu\})$ is the pre-image of $\mu$ by $L$ and $B_\M$ is the unit total variation ball.
We have $L^{-1}(B_\M)=L^+(B_\M)+K$ and we can identify $C_{K^\perp}$ with $L^+(B_\M)$. Since $L^+$ is bijective from $\M(\Omega)$ to $K^\perp$, the extreme points of $C_{K^\perp}$ are the image by $L^+$ of the Dirac masses.

The end of the proposition follows from Corollary \ref{cor:convex_fit} and from the fact that the lineality space of $\{u \in \B(\Omega), \|Lu\|_{\M(\Omega)}\leq 1\}$ is equal to $K$.
\end{proof}

Let us mention that, although the description of the extreme points follows directly from the results of Section 3, proving the existence of minimizers and the existence of extreme points is a considerably more difficult problem which needs a careful choice of topologies. The paper \cite{unser2017splines} provides a systematic way to construct Banach spaces and pseudo-inverse $L^+$ for ``spline admissible operators'' $L$ such as the fractional Laplacian. In addition, they prove existence of solutions by adding weak-* continuity assumptions on the sensing operator $\Phi$.}


\subsubsection{The total gradient variation}

Since its introduction in the field of image processing \cite{rudin1992nonlinear}, the total gradient variation proved to be an extremely valuable regularizer in diverse fields of data science and engineering. It is defined, for any locally integrable function $u$ as
\begin{equation*}
TV(u) \eqdef \sup\left(\int u\mathrm{div}(\phi) \, dx, \phi \in C^1_c(\RR^d)^d, \sup_{x\in \RR^d} \|\phi(x)\|_2\leq 1\right).
\end{equation*}
If the above quantity is finite, we say that $u$ has bounded variation and its gradient $D u$ is a Radon measure, with
\begin{equation*}
  TV(u) = \int_{\RR^d}|D u| = \|Du \|_{(\Mm(\RR^d))^d}.
\end{equation*}

Working in $\vecgal=L^{d/(d-1)}(\RR^d)$, one is led to consider the convex set $\cvx=\{ u\in \vecgal, TV(u)\leq 1\}$, referred to as the TV unit ball.

The generalized gradient operator is not a surjective operator. Hence, the analysis of \Cref{sec:prior} cannot help finding the extreme points of the TV ball. Still, those have been described in the fifties by Fleming in \cite{fleming1957functions} and refined analyses have been proposed more recently by Ambrosio, Caselles, Masnou and Morel in \cite{ambrosio2001connected}.

\begin{figure}
\begin{center}
\begin{tabular}{cc}
\includegraphics[height=4cm]{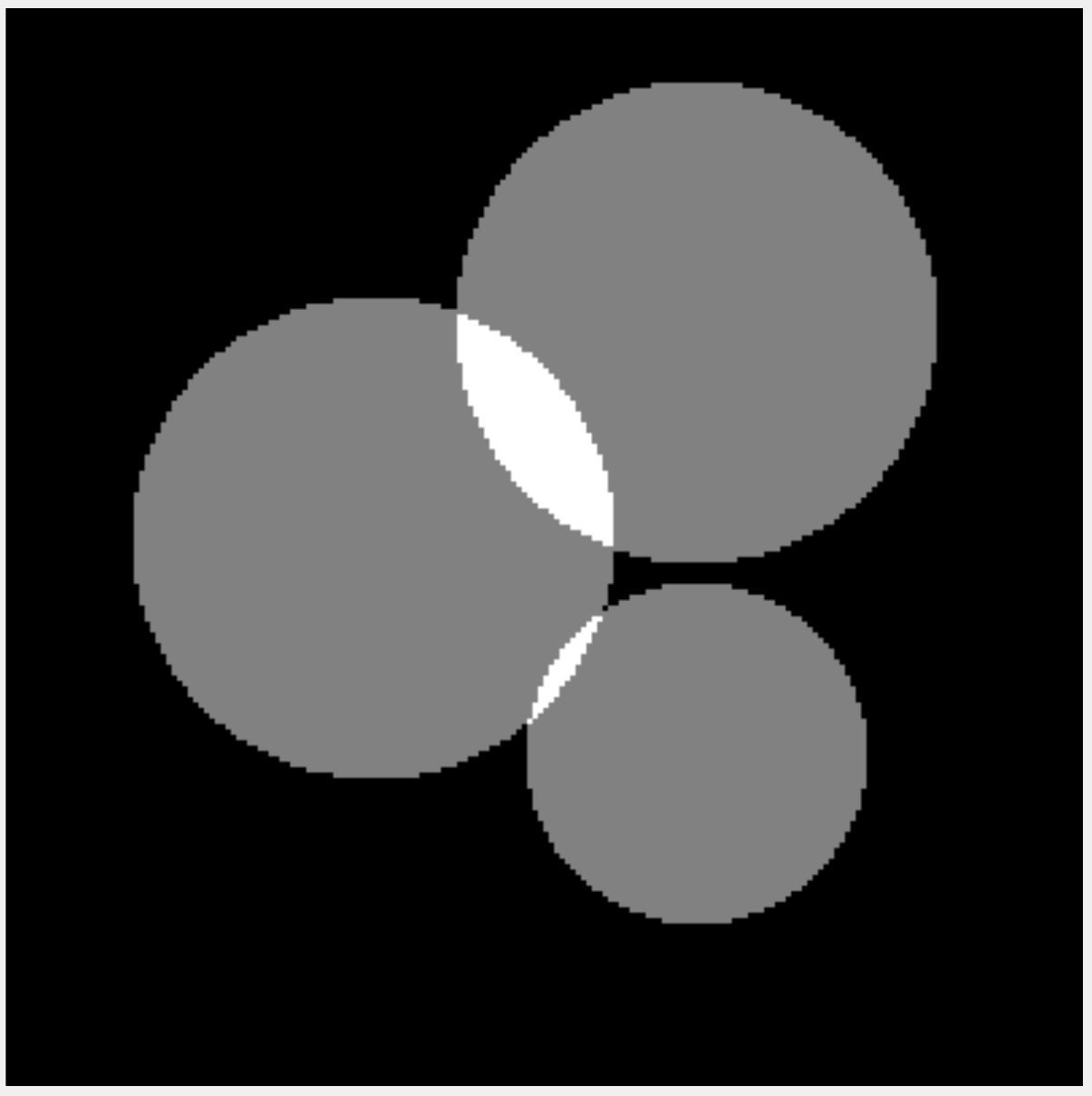} &
\includegraphics[height=4cm]{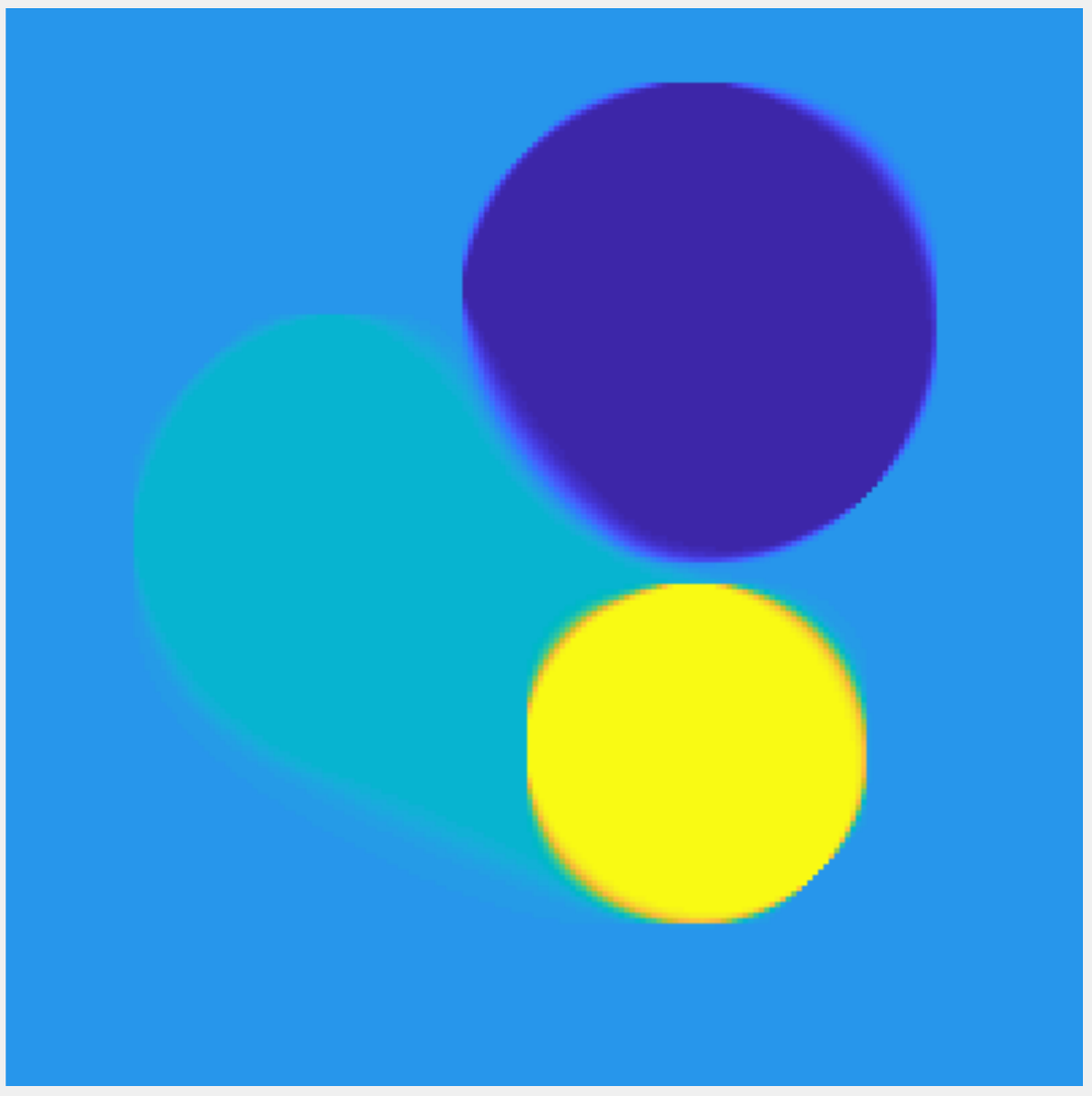} \\
{ (a)} & { (b)}
\end{tabular}
\caption{\label{fig:extremeTV}Illustration for the total gradient variation problem $\min \left\{ TV(u) : \Phi(u) =y \right\}$. Here, $\Phi$ is a linear mapping giving access to $3$ measurements, $y\in \R^3$, by performing the mean of an image $u$ of size $200 \times 200$ on $3$ different disks represented in (a). The TV problem is solved using a primal-dual algorithm, also known as the Chambolle-Pock algorithm \cite{Chambolle2011}. The recovered image is displayed in (b): it can be represented as the sum of $3$ indicator functions of simple sets. }
\end{center}
\end{figure}

\begin{theorem}[Extreme points of the TV ball \cite{fleming1957functions,ambrosio2001connected}]
\label{thm:extremeBV}
The extreme points of the unit TV unit ball are the indicators of simple sets normalized by their perimeter, i.e. functions of the form $u=\pm \frac{\one_F}{TV(\one_F)}$, where $F$ is an indecomposable and saturated subset of $\RR^d$.
\end{theorem}
Informally, the simple sets of $\RR^d$ are the simply connected sets with no hole. We refer the reader to \cite{ambrosio2001connected} for more details. Using \cref{thm:extremeBV} in conjunction with our results tell us that functions minimizing the total variation subject to a finite number of linear constraints can be expressed as a sum of a small number of indicators of simple sets, see for instance \cref{fig:extremeTV}, which is yet another theoretical result explaining the common observation that total variation tends to produce stair-casing \cite{nikolova2000local}.

\section{Conclusion}
In this paper we have developed representer theorems for convex regularized inverse problems \eqref{eq::mainproblem}, based on fundamental properties of the geometry of convex sets: the solution set can be entirely described using convex combinations of a small number of extreme points and extreme rays of the regularizer level set.

Obviously, the conclusion of Theorem~\ref{thm:first} is only nontrivial when $\minset$ has a ``sufficiently flat boundary'', in the sense that two or more faces of $\minset$ have dimension larger than $m$. For instance, if $\minset$ is strictly convex (\ie{} has only $0$-dimensional faces, except its interior\footnote{In this example, to simplify the discussion, we assume that $\vecgal$ has finite dimension.}), then the solution set $\sol$ is always reduced to a single extreme point of $\minset$! Nevertheless, several regularizers which are commonly used in the literature (notably sparsity-promoting ones) have that flatness property, and Theorem~\ref{thm:first} then provides interesting information on the structure of the solution set, as illustrated in Section~\ref{sec:app}.

To conclude, the structure theorem presented in this paper highlights the importance of describing the extreme points and extreme rays of the regularizer: this yields a fine description of the set of solutions of variational problems of the form~\eqref{eq::mainproblem}. Our theorem also suggests a principled way to design a regularizer. If a particular family of solutions is expected, then one may construct a suitable regularizer by taking the convex hull of this family. 
Finally, representer theorems have had a lot of success in the fields of approximation theory and machine learning \cite{scholkopf2001generalized}, in the frame of reproducible kernel Hilbert spaces. 
A reason of this success is that they allow to design efficient numerical procedures that yield \emph{infinite dimensional} solutions by solving \emph{finite dimensional} li\-near systems. 
Such numerical procedures have recently been extended to the case of Banach spaces for some simple instances of the problems described in this paper \cite{fernandez2016super,de2017exact,flinth2017exact}. The price to pay when going from a Hilbert space to a Banach space is that semi-infinite convex programs have to be solved instead of simpler linear systems.
We foresee that the results in this paper may help designing new efficient numerical procedures, since they allow to parameterize the solutions using extreme points and extreme rays only.

\section{Proofs of Section~\ref{sec:abstract}}
\label{sec:proof}

%
%

\subsection{Proof of Theorem~\ref{thm:first}}
\label{sec:prooffirst}
The set of solutions $\sol$ is precisely $\minset\cap \Phi^{-1}(\{y\})$, and the statement of the theorem amounts to describing its elementary faces.  Since $\Phi^{-1}(\{y\})$ is an affine space with codimension at most $m$, the main theorem of~\cite{klee_theorem_1963} almost provides the desired conclusion, but, for our particular case, it yields one extreme point/ray too many. Here is how to obtain the correct number. 

  Let $p$ be a point of $\sol$ such that $\face{p}{\sol}$ has dimension $j$. Up to a translation, it is not restrictive to assume that $p=0$, so that $\sol = \minset \cap \ker \Phi$.

  Let $\vecrest$ be the union of $\{0\}$ and all the lines $\ell$ such that $\minset\cap \ell$ contains an open interval which contains $0$. Note that $\vecrest$ is a linear space, the linear hull of $\face{0}{\minset}$.

We claim that $\codim_{\vecrest}\left(\vecrest\cap\ker \Phi\right)\leq m-1$. By contradiction assume that there is a complement $Z$ to $\vecrest\cap\ker\Phi$ in $\vecrest$ with dimension $m$. Then $\restr{\Phi}{Z}$ has rank $m$ and is a bijection, hence we may define
\begin{equation*}
  z\eqdef -\frac{\theta}{(1-\theta)}(\restr{\Phi}{Z})^{-1}(\Phi u_0) \in Z\subset \vecrest,
\end{equation*}
where $\theta\in ]0,1[$ and $u_0\in \minset$ is such that $\inf \reg \leq \reg(u_0)<\minval$. For $\theta$ small enough, $z\in \minset$, hence $\reg(z)\leq \minval$. Moreover, 
  \begin{equation}
    \Phi\left((1-\theta)z+\theta u_0 \right)= (1-\theta)\Phi z + \theta \Phi u_0 = 0.
  \end{equation}
 so that $(1-\theta)z+\theta u_0$ lies in {$\minset\cap \ker \Phi$}. Since $\reg(u_0) < \reg(z)$, and $\reg$ is convex
    \begin{equation}
    \label{eq:for_cvx_like}
      \reg\left((1-\theta)z+\theta u_0 \right)< \reg\left(z\right)\leq\minval, 
  \end{equation}
  we obtain a contradiction with the fact that $\minval$ is the minimal value of~\eqref{eq:thminreg}. As a result, $\codim_{\vecrest}\left(\vecrest\cap\ker \Phi\right)\leq m-1$.
  Observing that $\vecrest\cap\ker \Phi$ is the linear hull of $\face{0}{\sol}$, hence $j= \dim\left(\vecrest\cap\ker \Phi\right)$, we deduce that   
  \begin{equation}
    \dim \face{0}{\minset}\eqdef\dim \vecrest = \codim_{\vecrest}\left(\vecrest\cap\ker \Phi\right) + \dim\left(\vecrest\cap\ker \Phi\right) \leq m-1+j, 
  \end{equation}
and the first claim of the theorem is proved.

Now, applying the Carath\'eodory-Klee theorem (3) in~\cite{klee_theorem_1963}, $p$ is convex combination of at most $m+j$ (resp. $m+j-1$) extreme points (resp. extreme points or in an extreme ray) of $\face{0}{\minset}$.  The conclusion stems from the fact that the extreme points (resp. rays) of  $\face{0}{\minset}$ are extreme points (resp. rays) of $\minset$, see the proof of the main theorem in~\cite{klee_theorem_1963}.


\qed

\subsection{Proof of Corollary~\ref{coro:lines}}
\label{sec:proofconvex}
\begin{sloppypar}
  Now,  assume that the vector space ${\vecrec{}\eqdef\rec{\minset}\cap (-\rec{\minset})}$ is non trivial (otherwise the conclusion follows from Theorem~\ref{thm:first}). We note that for any $u\in \minset$, the convex function $v\mapsto \reg(u+v)$ is upper bounded by $\minval$ on $\vecrec{}$, hence is constant. As a result, possibly replacing $\reg$ with $\reg+\chi_{\minset}$, it is not restrictive to assume that $\reg$ is invariant by translation along $\vecrec{}$.
\end{sloppypar}

  Now, let $\projv$, $\projp$ be the canonical quotient maps and define $\tilde{\reg}$ and $\tilde{\Phi}$ by the commutative diagrams
  \[ \begin{tikzcd}
      \vecgal \arrow{r}{\reg} \arrow[swap]{d}{\projv} & \left[-\infty,\infty\right]  \\%
      \vecgal/\mvecrec \arrow[swap]{ru}{\tilde{\reg}}& 
\end{tikzcd}\qquad
 \begin{tikzcd}
\vecgal \arrow{r}{\Phi} \arrow[swap]{d}{\projv} & \RR^m  \arrow{d}{\projp} \\%
    \vecgal/\mvecrec \arrow{r}{\tilde{\Phi}}& \RR^m/\Phi(\mvecrec)
\end{tikzcd}\]
Note that $\tilde{\reg}$ is a convex function and that $\tilde{\Phi}$ is a linear map with rank $m-d$, where $d\eqdef \dim\left(\Phi(\mvecrec)\right)$. It is then natural to consider the problem 
\begin{equation}\label{eq:qthminreg}
    \min_{\tilde{u}\in \vecgal/\mvecrec}   \tilde{\reg}(\tilde{u})  \quad \mbox{s.t.}\quad \tilde{\Phi} \tilde{u}=\tilde{y},\tag{$\tilde{\Pp}$}
  \end{equation}
  where $\tilde{y}\eqdef\projp(y)$. In other words, one still wishes to minimize $\reg(u)$, but one is satisfied if the constraint $\Phi u=y$ merely holds up to an additional term $\Phi v$, where $v\in\mvecrec$. We observe that \eqref{eq:thminreg} and \eqref{eq:qthminreg} have the same value $\minval$, and the level set
  \begin{equation}
    \qminset\eqdef\enscond{\tilde{u}\in\vecgal/\mvecrec}{\tilde{\reg}(\tilde{u})\leq t^\star}=\projv(\minset)
  \end{equation}
  is convex linearly closed \emph{and contains no line}. Let $\qsol$ be the solution set to~\eqref{eq:qthminreg}. Theorem~\ref{thm:first} now describes the elements of the $j$-dimensional faces of $\qsol$ as convex combinations of $m-d+j$ (resp.\ $m-d+j-1$) extreme points (resp. extreme points or points in an extreme ray) of $\qminset$ that we denote by $\qe_1$,$\qe_2$,\ldots, $\qe_r$.

  To conclude, we have obtained $\projv(p)=\sum_i \theta_i \qe_i$, for some $\theta\in \RR_+^r$ with $\sum_i\theta_i=1$.
  Equivalently,  since 
$\qiso^{-1}(\cdot,0)$ provides one element in the corresponding class, this means that $p\in \qiso^{-1}(\sum_i \theta_i \qe_i,0)+\vecrec{}$. We get the claimed result by linearity of $\qiso^{-1}$.\qed


\begin{remark}
  Incidentally, we note that for $\cible=\{y\}$, the face $\face{p}{\sol_{\{y\}}}$ is isomorphic (through $\qiso$) to $\face{\projv(p)}{\projv(\sol_{\{y\}})}\times (\vecrec{}\cap \ker\Phi)$.
\end{remark}

\section*{Acknowledgments}
This work was initially started by two different groups composed of A. Chambolle, F. de Gournay and P. Weiss
 on one side and C. Boyer, Y. De Castro and V. Duval on the other side. The authors realized that they were working on a similar topic when a few of them met at the Cambridge semester in mathematical imaging during the Isaac Newton Institute (Cambridge) semester ``Variational methods and effective algorithms for imaging and vision'', supported by EPSRC Grant N.~EP/K032208/1. They therefore decided to write a joint paper and wish to thank the organizers for giving them this opportunity.
 This work initially started with the help of T. Pock through a few numerical experiments and with discussions with J. Fadili and C. Poon. 
 
 The work of A.C. was partially supported by a grant of the Simons Foundation.
\bibliographystyle{abbrv}
\bibliography{Biblio}

\begin{thebibliography}{10}

\bibitem{adcock2016generalized}
B.~Adcock and A.~C. Hansen.
\newblock Generalized sampling and infinite-dimensional compressed sensing.
\newblock {\em Foundations of Computational Mathematics}, 16(5):1263--1323,
  2016.

\bibitem{ahmed2014blind}
A.~Ahmed, B.~Recht, and J.~Romberg.
\newblock Blind deconvolution using convex programming.
\newblock {\em IEEE Transactions on Information Theory}, 60(3):1711--1732,
  2014.

\bibitem{aliprantis_infinite_2006}
C.~Aliprantis and K.~Border.
\newblock {\em Infinite {Dimensional} {Analysis}: {A} {Hitchhiker}'s {Guide}}.
\newblock Springer Berlin Heidelberg, 2006.

\bibitem{ambrosio2001connected}
L.~Ambrosio, V.~Caselles, S.~Masnou, and J.-M. Morel.
\newblock Connected components of sets of finite perimeter and applications to
  image processing.
\newblock {\em Journal of the European Mathematical Society}, 3(1):39--92,
  2001.

\bibitem{barvinok_problems_1995}
A.~I. Barvinok.
\newblock Problems of distance geometry and convex properties of quadratic
  maps.
\newblock {\em Discrete \& Computational Geometry}, 13(2):189--202, Mar. 1995.

\bibitem{Beurling1938}
A.~Beurling.
\newblock Sur les int{\'e}grales de {F}ourier absolument convergentes et leur
  application {\`a} une transformation fonctionnelle.
\newblock In {\em Ninth Scandinavian Mathematical Congress}, pages 345--366,
  1938.

\bibitem{beutler}
F.~J. Beutler.
\newblock The operator theory of the pseudo-inverse i. bounded operators.
\newblock 10, 06 1965.

\bibitem{bourbaki_espaces_2007}
N.~Bourbaki.
\newblock {\em Espaces vectoriels topologiques: {Chapitres} 1{\`a} 5}.
\newblock Springer-Verlag, Berlin Heidelberg, 2007.

\bibitem{bredies_sparsity_2018}
K.~Bredies and M.~Carioni.
\newblock Sparsity of solutions for variational inverse problems with
  finite-dimensional data.
\newblock {\em arXiv:1809.05045 [math]}, Sept. 2018.
\newblock arXiv: 1809.05045.

\bibitem{candes2014towards}
E.~J. Cand{\`e}s and C.~Fernandez-Granda.
\newblock Towards a mathematical theory of super-resolution.
\newblock {\em Communications on Pure and Applied Mathematics}, 67(6):906--956,
  2014.

\bibitem{candes2009exact}
E.~J. Cand{\`e}s and B.~Recht.
\newblock Exact matrix completion via convex optimization.
\newblock {\em Foundations of Computational mathematics}, 9(6):717, 2009.

\bibitem{Chambolle2011}
A.~Chambolle and T.~Pock.
\newblock A first-order primal-dual algorithm for convex problems with
  applications to imaging.
\newblock {\em Journal of Mathematical Imaging and Vision}, 40(1):120--145, May
  2011.

\bibitem{chandrasekaran2012convex}
V.~Chandrasekaran, B.~Recht, P.~A. Parrilo, and A.~S. Willsky.
\newblock The convex geometry of linear inverse problems.
\newblock {\em Foundations of Computational mathematics}, 12(6):805--849, 2012.

\bibitem{daniilidis2007some}
A.~Daniilidis and Y.~G. Ramos.
\newblock Some remarks on the class of continuous (semi-) strictly quasiconvex
  functions.
\newblock {\em Journal of optimization theory and applications}, 133(1):37--48,
  2007.

\bibitem{dattorro_convex_2005}
J.~Dattorro.
\newblock {\em Convex {Optimization} \& {Euclidean} {Distance} {Geometry}}.
\newblock Meboo Publishing, 2005.

\bibitem{de2017exact}
Y.~De~Castro, F.~Gamboa, D.~Henrion, and J.-B. Lasserre.
\newblock Exact solutions to super resolution on semi-algebraic domains in
  higher dimensions.
\newblock {\em IEEE Transactions on Information Theory}, 63(1):621--630, 2017.

\bibitem{Donoho2006Compressed}
D.~L. Donoho.
\newblock Compressed sensing.
\newblock {\em IEEE Transactions on Information Theory}, 52(4):1289--1306,
  2006.

\bibitem{donoho2005sparse}
D.~L. Donoho and J.~Tanner.
\newblock Sparse nonnegative solution of underdetermined linear equations by
  linear programming.
\newblock {\em Proceedings of the National Academy of Sciences of the United
  States of America}, 102(27):9446--9451, 2005.

\bibitem{drusvyatskiy2015extreme}
D.~Drusvyatskiy, S.~A. Vavasis, and H.~Wolkowicz.
\newblock Extreme point inequalities and geometry of the rank sparsity ball.
\newblock {\em Mathematical Programming}, 152(1-2):521--544, 2015.

\bibitem{dubins1962extreme}
L.~E. Dubins.
\newblock On extreme points of convex sets.
\newblock {\em Journal of Mathematical Analysis and Applications},
  5(2):237--244, 1962.

\bibitem{Duval2015Exact}
V.~Duval and G.~Peyr\'e.
\newblock Exact support recovery for sparse spikes deconvolution.
\newblock {\em Foundations of Computational Mathematics}, 15(5):1315--1355,
  2015.

\bibitem{fernandez2016super}
C.~Fernandez-Granda.
\newblock Super-resolution of point sources via convex programming.
\newblock {\em Information and Inference: A Journal of the IMA}, 5(3):251--303,
  2016.

\bibitem{fisher_spline_1975}
S.~D. Fisher and J.~W. Jerome.
\newblock Spline solutions to {L}1 extremal problems in one and several
  variables.
\newblock {\em Journal of Approximation Theory}, 13(1):73--83, 1975.

\bibitem{fleming1957functions}
W.~Fleming.
\newblock Functions with generalized gradient and generalized surfaces.
\newblock {\em Annali di Matematica Pura ed Applicata}, 44(1):93--103, 1957.

\bibitem{flinth2017exact}
A.~Flinth and P.~Weiss.
\newblock Exact solutions of infinite dimensional total-variation regularized
  problems.
\newblock {\em arXiv preprint arXiv:1708.02157}, 2017.

\bibitem{Fogel2015Convex}
F.~Fogel, R.~Jenatton, F.~Bach, and A.~D'Aspremont.
\newblock Convex relaxations for permutation problems.
\newblock {\em Advances in Neural Information Processing Systems},
  36(4):1016--1024, 2015.

\bibitem{gupta2018continuous}
H.~Gupta, J.~Fageot, and M.~Unser.
\newblock Continuous-domain solutions of linear inverse problems with tikhonov
  vs. generalized tv regularization.
\newblock {\em arXiv preprint arXiv:1802.01344}, 2018.

\bibitem{hiriart-urruty_convex_1993}
J.-B. Hiriart-Urruty and C.~Lemar{\'e}chal.
\newblock {\em Convex {Analysis} and {Minimization} {Algorithms} {I}}, volume
  305 of {\em Grundlehren der mathematischen {Wissenschaften}}.
\newblock Springer Berlin Heidelberg, Berlin, Heidelberg, 1993.

\bibitem{klee_theorem_1963}
V.~Klee.
\newblock On a theorem of {Dubins}.
\newblock {\em Journal of Mathematical Analysis and Applications},
  7(3):425--427, Dec. 1963.

\bibitem{klee_extremal_1957}
V.~L. Klee.
\newblock Extremal structure of convex sets.
\newblock {\em Archiv der Mathematik}, 8(3):234--240, Aug. 1957.

\bibitem{Krein1938}
M.~Kre{\v{i}}n.
\newblock The {L}-problem in an abstract normed linear space.
\newblock In I.~Ahiezer and M.~Kre{\v{i}}n, editors, {\em Some questions in the
  theory of moments}, chapter~4. Gos. Nau{\v{c}}no-Tehn. Izdat. Ukraine, 1938.
\newblock English Transl. Amer. Math. Soc., Providence, R.I., 1962. MR 29 {\#}
  5073.

\bibitem{krein_markov_1977}
M.~G. Krein and A.~A. Nudelman.
\newblock {\em The {Markov} moment problem and extremal problems: ideas and
  problems of {P}. {L}. {Cebysev} and {A}. {A}. {Markov} and their further
  development}.
\newblock Number v. 50 in Translations of mathematical monographs. American
  Mathematical Society, Providence, R.I, 1977.

\bibitem{lawson1995solving}
C.~L. Lawson and R.~J. Hanson.
\newblock {\em Solving least squares problems}, volume~15.
\newblock Siam, 1995.

\bibitem{matousek2007understanding}
J.~Matousek and B.~G{\"a}rtner.
\newblock {\em Understanding and using linear programming}.
\newblock Springer Science \& Business Media, 2007.

\bibitem{mcmullen1970maximum}
P.~McMullen.
\newblock The maximum numbers of faces of a convex polytope.
\newblock {\em Mathematika}, 17(2):179--184, 1970.

\bibitem{motzkin1957comonotone}
T.~S. Motzkin.
\newblock Comonotone curves and polyhedra.
\newblock {\em Bull. Amer. Math. Soc.}, 63, 1957.

\bibitem{nikolova2000local}
M.~Nikolova.
\newblock Local strong homogeneity of a regularized estimator.
\newblock {\em SIAM Journal on Applied Mathematics}, 61(2):633--658, 2000.

\bibitem{pataki2000geometry}
G.~Pataki.
\newblock The geometry of semidefinite programming.
\newblock In {\em Handbook of semidefinite programming}, pages 29--65.
  Springer, 2000.

\bibitem{Richard2013Intersecting}
E.~Richard, F.~Bach, and J.~P. Vert.
\newblock Intersecting singularities for multi-structured estimation.
\newblock In {\em International Conference on International Conference on
  Machine Learning}, pages III--1157, 2013.

\bibitem{richard2014tight}
E.~Richard, G.~R. Obozinski, and J.-P. Vert.
\newblock Tight convex relaxations for sparse matrix factorization.
\newblock In Z.~Ghahramani, M.~Welling, C.~Cortes, N.~D. Lawrence, and K.~Q.
  Weinberger, editors, {\em Advances in Neural Information Processing Systems
  27}, pages 3284--3292. Curran Associates, Inc., 2014.

\bibitem{rudin1992nonlinear}
L.~I. Rudin, S.~Osher, and E.~Fatemi.
\newblock Nonlinear total variation based noise removal algorithms.
\newblock {\em Physica D: nonlinear phenomena}, 60(1-4):259--268, 1992.

\bibitem{scholkopf2001generalized}
B.~Sch{\"o}lkopf, R.~Herbrich, and A.~J. Smola.
\newblock A generalized representer theorem.
\newblock In {\em International conference on computational learning theory},
  pages 416--426. Springer, 2001.

\bibitem{Scholkopf2002Learning}
B.~Scholkopf and A.~J. Smola.
\newblock {\em Learning with kernels :}.
\newblock MIT Press,, 2002.

\bibitem{shapiro2001duality}
A.~Shapiro.
\newblock On duality theory of conic linear problems.
\newblock In {\em Semi-infinite programming}, pages 135--165. Springer, 2001.

\bibitem{shohat1943problem}
J.~A. Shohat and J.~D. Tamarkin.
\newblock {\em The problem of moments}.
\newblock Number~1. American Mathematical Soc., 1943.

\bibitem{Tang2013Compressed}
G.~Tang, B.~N. Bhaskar, P.~Shah, and B.~Recht.
\newblock Compressed sensing off the grid.
\newblock {\em IEEE Transactions on Information Theory}, 59(11):7465--7490,
  2013.

\bibitem{unser2016representer}
M.~Unser, J.~Fageot, and H.~Gupta.
\newblock {Representer Theorems for Sparsity-Promoting $\ell_{1}$
  Regularization}.
\newblock {\em IEEE Transactions on Information Theory}, 62(9):5167--5180,
  2016.

\bibitem{unser2017splines}
M.~Unser, J.~Fageot, and J.~P. Ward.
\newblock Splines are universal solutions of linear inverse problems with
  generalized tv regularization.
\newblock {\em SIAM Review}, 59(4):769--793, 2017.

\bibitem{Wendland2005Scattered}
H.~Wendland.
\newblock {\em Scattered data approximation}.
\newblock Cambridge University Press, 2005.

\bibitem{Zuhovickii1948}
S.~Zuhovicki\u{i}.
\newblock Remarks on problems in approximation theory.
\newblock {\em Mat. Zbirnik KDU}, pages 169--183, 1948.
\newblock (Ukrainian).

\bibitem{Zuhovickii1962}
S.~Zuhovicki\u{i}.
\newblock On approximation of real functions in the sense of {P.L.}
  \v{C}eby\v{s}ev.
\newblock {\em AMS Translations of Mathematical Monographs.}, 19(2):221--252,
  1962.

\end{thebibliography}

\end{document}